\documentclass[a4paper,pdftex,reqno]{amsart}

\usepackage{geometry}
\title{Separability criteria for loops via the Goldman bracket}
\author{Aoi Wakuda}

\makeatletter
\@namedef{subjclassname@2020}{\textup{2020} Mathematics Subject Classification}
\makeatother

\address{GRADUATE SCHOOL OF MATHEMATICAL SCIENCES, UNIVERSITY OF TOKYO, 3-8-1 KOMABA, MEGURO-KU, TOKYO, 153-8914,
JAPAN}
\email{aoichan19991226@g.ecc.u-tokyo.ac.jp}

\usepackage{amsmath} 
\usepackage{amssymb} 
\usepackage{amsfonts} 
\usepackage{mathrsfs} 
\usepackage{amsthm} 
\usepackage{bm} 
\usepackage[new]{old-arrows} 
\usepackage[dvipdfmx]{graphicx} 
\usepackage{wrapfig} 
\usepackage{color} 
\usepackage[labelsep=colon]{caption} 
\usepackage[labelformat=empty,subrefformat=parens]{subcaption} 
\usepackage{multicol} 
\usepackage{units} 
\usepackage[all]{xy} 
\usepackage{tikz-cd} 
\usepackage{hyperref} 
\usepackage{color} 
\usepackage{latexsym}
\usepackage{mathtools}
\usepackage{comment}
\usepackage{tikz-cd}
\usepackage{tikz}
\usepackage{float}
\usepackage{enumerate}
\usetikzlibrary{decorations.markings}
\usetikzlibrary{arrows,arrows.meta}
\usepackage{graphicx}
\usepackage{overpic}

\newtheorem{thm}{Theorem}[section]

\newtheorem{prop}[thm]{Proposition}
\newtheorem{lemma}[thm]{Lemma}
\newtheorem{cor}[thm]{Corollary}

\newtheorem*{mthm*}{Main Theorem}

\theoremstyle{definition}

\newtheorem{ex}[thm]{Example}
\newtheorem{remark}[thm]{Remark}

\newcommand{\C}{\mathbb{C}}

\newcommand{\Al}{\alpha}
\newcommand{\Be}{\beta}

\begin{document}

\begin{abstract}
We provide some explicit algebraic criteria in terms of the Goldman bracket to decide whether two free homotopy classes of loops on an oriented surface admit disjoint representatives. We extend the method of Kabiraj \cite{Kabiraj2016} using the hyperbolic geometry of surfaces to prove these criteria. As an application, we show that the center of the Goldman Lie algebra of a pair of pants is generated by the class of the constant loop together with the classes of loops that wind multiple times around a single puncture or boundary component. This case was not covered by Kabiraj \cite{Kabiraj2016}, since a pair of pants is not filled by simple closed curves.
\end{abstract}
\maketitle

\section{Introduction}
Let \( \Sigma \) be a connected oriented surface (possibly with boundary and punctures). We assume the Euler characteristic of \( \Sigma \) is negative so that \( \Sigma \) admits a complete hyperbolic metric. Denote by \( \hat{\pi} \) the set of free homotopy classes of loops on \( \Sigma \). Let \( K \) be a unital commutative ring of characteristic zero. For any set \( X \), we denote by \( KX \) the free \( K \)-module generated by \( X \).

In the 1980s, Goldman \cite{Goldman1986} defined a Lie bracket on \( K\hat{\pi} \) based on the intersection points of generic representatives of two free homotopy classes. This Lie algebra is now known as the \emph{Goldman Lie algebra}, and its bracket is called the \emph{Goldman bracket}. This bracket reflects the intersection points of two free loops on the surface. For $x, y \in \hat{\pi}$, the geometric intersection number $i(x, y)$ is defined as the minimal number of transversal double points among all pairs of generic representatives of $x$ and $y$. That is,
\[
i(x, y) := \min\{ \#(a \cap b) \mid a \in x, b \in y, \text{ and $a,b$ are generic} \}.
\]

In particular, \( i(x, y) = 0 \) means \( x \) and \( y \) admit disjoint representatives, which we call the separability of $x$ and $y$.

The surface \( \Sigma \) admits a complete hyperbolic metric. For simple closed curves, geodesic representatives realize the geometric intersection number. For non-simple loops, the same holds after refining the notion of intersection points (see Subsection~\ref{A generalized notion of intersection points} for details). By using geodesic representatives of simple closed curves that fill the surface, Kabiraj~\cite{Kabiraj2016} determined the center of the Goldman Lie algebra for surfaces with boundary, except for a pair of pants. As shown in this paper, this geometric approach can be extended to loops that are not necessarily simple. As an application, we obtain several theorems (Theorems~\ref{WSC1}, \ref{SSC1}, and~\ref{COPOP1}). In particular, we determine the center of the Goldman Lie algebra of a pair of pants in Theorem~\ref{COPOP1}.

To study the separability of two free homotopy classes of loops, we begin with the following theorem by Goldman \cite{Goldman1986}.

\begin{thm}\cite[Theorem 5.17]{Goldman1986}
\label{Goldman1986_1}
Let \( x, y \in \hat{\pi} \), where \( x \) is represented by a simple closed curve. Then \( [x, y] = 0 \) in \( \mathbb{Z}\hat{\pi} \) if and only if \( i(x, y) = 0 \).
\end{thm}

The simplicity of $x \in \hat\pi$ is necessary in Theorem 5.17 of \cite{Goldman1986}.
In general, the vanishing of the Goldman bracket does not imply that of the geometric intersection number. In particular, for $x,y \in \hat{\pi}$, the bracket $[x,y]=0$ holds trivially when $x=y$. We call it the \emph{trivial vanishing case of the bracket}: the bracket $[x,x]$ vanishes by antisymmetry, while the geometric intersection number $i(x,x)$ may be nonzero. Note that some classes with self-intersection points satisfy $i(x,x) > 0$. Moreover, Chas \cite{Chas2010} gave some explicit $x,y \in\hat \pi$ satisfying $[x, y] = 0, x\neq y$ and $i(x, y) > 0$. See Example 9.1 of \cite{Chas2010} for details.

These observations motivate our search for an algebraic characterization of the separability of (not necessarily simple) loops in terms of the Goldman bracket. To state our first result, we introduce the following notation and definition. For each \( n \in \mathbb{N}_{\geq 1} \), let \( f_n : S^1 \to S^1 \) be the map defined by \( f_n(z) = z^n \). For \( m \in \mathbb{Z} \setminus \{0\} \), we denote by \( x^m \in \hat{\pi} \) the free homotopy class represented by the composition \( \gamma \circ f_m \), where \( \gamma : S^1 \to \Sigma \) is a representative of \( x \). Our first theorem provides such a criterion.

\begin{thm}[Theorem \ref{WSC2}] \label{WSC1}
Let \( x, y \in \hat{\pi} \) and \( m \in \mathbb{N}_{\geq 2} \).  
Then \( [x^m, y] = 0 \) in \( K\hat{\pi} \) if and only if \( i(x, y) = 0 \) or \( y = x^m \).
\end{thm}

The Goldman bracket detects also information about self-intersection of loops on surfaces. An element \( x \in \hat{\pi} \) is \emph{primitive} if any representative of \( x \) does not factor through \( f_n \) for any \( n \in \mathbb{N}_{\geq 1} \). Motivated by computer experiments, Chas and Kabiraj \cite{Chas-Kabiraj2023} gave the following conjecture: if $x \in \hat\pi$ is a primitive class \( x \), then the number of terms (counted with multiplicity) appearing in the bracket \( [x, x^n] \), for \( n \geq 2 \) or \( n = -1 \), is equal to \( |2n| \) times the self-intersection number of \( x \). The conjecture was partially solved in \cite{Chas-Krongold2010} and \cite{Chas-Kabiraj2023}. Chas and Krongold \cite{Chas-Krongold2010} proved it in the case of surfaces with boundary and \( n \geq 3 \). Moreover, Chas and Kabiraj~\cite{Chas-Kabiraj2023} proved that a primitive class $x$ contains a simple representative if and only if one of the following holds:
\begin{enumerate}
    \item \( [x, x^n] = 0 \) for some \( n \in \{2, 3, \dots\} \),
    \item \( [x, \bar{x}] = 0 \), where \( \bar{x} \) denotes the conjugacy class of \( x \) with the opposite orientation.
\end{enumerate}
Our work refines and extends the methods of \cite{Chas-Kabiraj2023}, and as a consequence, condition (1) above arises naturally as a corollary of Theorem~\ref{WSC1} by substituting $x$ for $y$ in Theorem~\ref{WSC1}. As for condition (2), it was recently proved in \cite{Alonso-Paternain-Peraza-Reisenberger2023}, using a generalization of the combinatorial techniques developed by Chas in \cite{Chas2004}.

We now return to the topic of separability. Theorem \ref{WSC1} does not exclude the trivial vanishing case of the bracket (\( y = x^m \)), but does our second theorem, which provides complete algebraic criteria for the separability of loops.

\begin{thm}[Theorem \ref{SSC2}] \label{SSC1}
Let $x, y \in \hat{\pi}$. The following four conditions are equivalent,\\
{\rm (1)} $i(x, y) = 0$,\\
{\rm (2)} There exist distinct $m_1, m_2 \in \mathbb{N}_{\geq 1}$ such that $[x^{m_1}, y] = [x^{m_2}, y] = 0 $ in \( K\hat{\pi} \),\\
{\rm (3)} There exist distinct $m_1, m_2 \in \mathbb{N}_{\geq 1}$ such that $[x^{m_1}, y] = [x, y^{m_2}] = 0$ in \( K\hat{\pi} \),\\
{\rm (4)} There exists $m \in \mathbb{N}_{\geq 2}$ and non-zero divisors $c_1, c_2 \in K$ such that $[x^m, c_1 y + c_2 y^{-1}] = 0$ in \( K\hat{\pi} \).
\end{thm}
\noindent Here we remark $[x^m, x^n]$ does not vanish in general if $m \neq n$ and $m, n \neq 0$.

As for a surface $\Sigma$ whose Euler characteristic is nonnegative, then \( \Sigma \) is either the sphere, disk, annulus or torus. In the first three cases, the Goldman Lie algebra is trivial.  For the torus, the vanishing of the Goldman bracket of two free homotopy classes implies that of their geometric intersection number. For a proof in the torus case, see Lemma~7.6 of \cite{Chas2010}.

The proof of Theorem \ref{SSC1} is applied to compute the center of the Goldman Lie algebra. We call a loop on \( \Sigma \) \emph{non-essential} if it is the constant loop or is freely homotopic to a single puncture or a boundary component, and \emph{essential} otherwise. Studying the moduli spaces of flat bundles, Etingof \cite{Etingof2006} proved that the center of the Goldman Lie algebra $\C\hat{\pi}$ on a closed surface is generated by the class of the constant loop. Kawazumi and Kuno \cite{Kawazumi-Kuno2013} used representation theory of the symplectic group to prove that the center of the Goldman Lie algebra $\mathbb{Q}\hat{\pi}$ on a surface of infinite genus and one boundary component is generated by the class of non-essential loops. In \cite{Kabiraj2016}, Kabiraj proved that the center of the Goldman Lie algebra on surfaces with boundary, except a pair of pants, is generated by the class of non-essential loops. His proof is based on hyperbolic geometry and relies on the assumption that the surface can be filled by simple closed geodesics with respect to a complete hyperbolic metric, so it covers all surfaces with negative Euler characteristic other than a pair of pants. In this paper, we extend his approach to handle closed geodesics that are not necessarily simple and obtain the following result.

\begin{thm}[Theorem \ref{COPOP2}] \label{COPOP1}
The center of the Goldman Lie algebra of a pair of pants is generated by the class of non-essential loops as a $K$-module.
\end{thm}

In the body of this paper, these theorems will be proved as follows. From the definition of the Goldman bracket, if \( i(x,y) = 0 \) or \( y = x^{m} \), then \([x^{m}, y] = 0\). Theorem~\ref{WSC1} (Theorem \ref{WSC2}) means the converse direction. In order to prove it, we fix a complete hyperbolic metric \( X \) on the surface \( \Sigma \), and take the geodesic representatives \( x^{m}(X) \) and \( y(X) \) for \( x^{m} \) and \( y \in \hat{\pi} \), respectively. We may assume the geodesics intersect transversely. Now we focus on an intersection point \( P \) whose forward angle \( \phi_{P} \) is minimum among all the intersection points \( x^{m}(X) \cap y(X) \). Since \([x^{m}, y] = 0\), the term coming from \( P \) must be cancelled. We take another intersection point \( Q \) whose term in \([x^{m}, y]\) is cancelled with that of \( P \). Next, we construct a bi-infinite zigzag curve \( C \) on the hyperbolic plane by first taking a lift of \( P \) to the hyperbolic plane and then alternately lifting \( x^{m}(X) \) and \( y(X) \) from that point. Similarly, for the point \( Q \), we construct another bi-infinite zigzag curve \( D \) on the hyperbolic plane. The zigzag curves \( C \) and \( D \) may have various relative positions to each other. For any relative positions, we will find an intersection point \( R \in x^{m}(X) \cap y(X) \) whose forward angle $\phi_{R}$ is strictly smaller than \( \phi_{P} \). This contradicts the choice of the point \( P \). Hence \( x^{m}(X) \) and \( y(X) \) have no intersection points, and Theorem~\ref{WSC1} (Theorem \ref{WSC2}) follows.

It is a new contribution of this paper to study the case where a lift of $x^{m}(X)$ in $C$ and a lift of $x^{m}(X)$ in $D$ intersect transversely. In this case the closed geodesic $x^{m}(X)$ has a self-intersection point, and the assumption $m \ge 2$ guarantees the existence of an intersection point $R \in x^{m}(X)\cap y(X)$ whose forward angle $\phi_{R}$ is strictly smaller than $\phi_{P}$. (See Case~(I) of Lemma~\ref{xyz_-} for details.) Theorem~4.5 of \cite{Chas-Kabiraj2023} already studied the special case $y=x$, which also appears as Corollary~\ref{1_m} in this paper.

Now we look at a basic fact that $x^{m}=x^{n}$ in $\hat{\pi}$ for a nontrivial loop $x$ implies $m=n$. We combine the fact with a similar argument stated above to prove Theorem~\ref{SSC1} (Theorem \ref{SSC2}).

Since non-essential loops admit representatives disjoint from all other loops, linear combinations of non-essential loops are in the center.
The converse direction is an essential part of Theorem \ref{COPOP1} (Theorem \ref{COPOP2}) For the converse direction, consider a linear combination $y=\sum_{j=1}^{k} c_{j} y_{j}$ with all coefficients $c_{j}$ nonzero. Let $x$ be a figure-eight curve on a pair of pants, and assume $I \coloneqq i(x,y_{1}) + \cdots + i(x,y_{k}) \ge 1$. We may also assume $i(x,y_{1}) \ge 1$. Since $y$ is in the center, we have $[x^{m},y]=0$ for every natural number $m \ge 2$. As before, we focus on an intersection point $P \in x^{m}(X)\cap y_1(X)$ whose forward angle $\phi_{P}$ is minimum among all the intersection points \( x^{m}(X) \cap y_1(X) \). By a similar argument stated above, we obtain an intersection point $R \in x^{m}(X)\cap y_1(X)$ whose forward angle $\phi_{R}$ is strictly smaller than $\phi_{P}$, which is impossible. Therefore we obtain $i(x,y_{j})=0$ for all $j$, and hence each $y_{j}$ is non-essential.

\noindent
\textbf{Acknowledgment.} The author is deeply grateful to his supervisor, Nariya Kawazumi, for many valuable discussions. The author also wishes to thank Kento Sakai and Toyo Taniguchi for numerous inspiring conversations. This research was supported by the WINGS-FMSP program of the Graduate School of Mathematical Sciences, University of Tokyo.

\vspace{1em}
\noindent\textbf{\Large Contents}
\vspace{1em}

\noindent\textbf{1\quad Introduction} \dotfill \textbf{1} \\[0.3em]

\noindent\textbf{2\quad Preliminaries on the Goldman bracket and hyperbolic geometry} \dotfill \textbf{3} \\[0.3em]
\hspace*{1.8em}2.1\quad The definition of the Goldman Lie algebra \dotfill 4 \\[0.3em]
\hspace*{1.8em}2.2\quad A generalized notion of intersection points \dotfill 4 \\[0.3em]
\hspace*{1.8em}2.3\quad Closed geodesics on a complete hyperbolic surface \dotfill 4 \\[0.3em]
\hspace*{1.8em}2.4\quad Hyperbolic isometries \dotfill 5 \\[0.3em]

\noindent\textbf{3\quad Separability Criteria} \dotfill \textbf{6} \\[0.3em]
\hspace*{1.8em}3.1\quad A zigzag curve induced by an intersection point of two closed geodesics \dotfill 6 \\[0.3em]
\hspace*{1.8em}3.2\quad Intersection of two zigzag curves \dotfill 6 \\[0.3em]
\hspace*{1.8em}3.3\quad Key results \dotfill 11 \\[0.3em]
\hspace*{1.8em}3.4\quad Proof of separability criteria \dotfill 15 \\[0.3em]

\noindent\textbf{4\quad The center of the Goldman Lie algebra of a pair of pants} \dotfill \textbf{18} \\[0.3em]

\section{Preliminaries on the Goldman bracket and hyperbolic geometry} \label{Goldman Lie algebra and hyperbolic geometry}

In this section, we recall some facts on the Goldman bracket and hyperbolic geometry which will be used for the proof of the separability criteria.

\subsection{The definition of the Goldman Lie algebra} \label{The Goldman Lie algebra}

Let \( \Sigma \) be a connected oriented surface and $K$ a unital commutative ring of characteristic zero. Denote by $\hat{\pi}$ the set of free homotopy classes of free loops on \( \Sigma \) and by $|{x}| \in \hat\pi$ the free homotopy class of a loop $x$.

The {\it Goldman bracket} of \(x, y \in \hat{\pi}\) is defined by
\begin{align}
[x,y] \coloneqq \sum_{P\in{x}\cap{y}}\varepsilon_{P}(x,y)|x_{P}y_{P}|.
\end{align}
Here the representatives $x$ and $y$ are chosen so that they intersect transversely in a set of double points ${x}\cap{y}$, $\varepsilon_{P}(x,y)$ denotes the sign of the intersection between $x$ and $y$ at an intersection point ${P}$, and ${x_{P}y_{P}}$ denotes the loop product of $x$ and $y$ at $P$. In \cite{Goldman1986}, Goldman proved the bracket defined above is well-defined, skew-symmetric, and satisfies the Jacobi identity on $K\hat{\pi}$. In other words, $K\hat{\pi}$ is a Lie algebra.

\subsection{A generalized notion of intersection points} \label{A generalized notion of intersection points}

In Sections~3 and 4, we compute the Goldman bracket of two free homotopy classes of loops by taking their geodesic representatives with respect to a fixed complete hyperbolic metric on a surface. However, geodesic representatives do not always intersect only at double points, as required in Goldman's original definition. To compute the Goldman bracket using geodesic representatives, Kabiraj~\cite{Kabiraj2018} first introduced a notion of intersection points that allows us to handle cases where several intersection points occur at the same point on the surface. 
Later, Chas–Kabiraj~\cite{Chas-Kabiraj2022} gave a simpler and more convenient formulation of this idea which we will adopt through this paper. 

If $a$ and $b$ are two curves intersecting transversely, we call a point $P \in a \cap b$, together with a choice of a pair of small arcs, one from $a$ and the other from $b$ intersecting only at $P$, an $(a,b)$-\textit{transverse intersection point} (in the sense of~\cite{Chas-Kabiraj2022}). The original definition of the Goldman bracket can then be extended to allow summation over all $(a, b)$-transverse intersection points. For further details, see~\cite[Appendix~B]{Kabiraj2018} and~\cite[Section~2]{Chas-Kabiraj2022}.

Moreover, the geometric intersection number can be reformulated in terms of $(a,b)$-\textit{transverse intersection points}. For any two free homotopy classes $x, y \in \hat{\pi}$, we have
\[
i(x, y) = \min \{ \#(a, b)\text{-transverse intersection points} \mid 
a \in x,\; b \in y,\; a \text{ and } b \text{ intersect transversely} \}.
\]
Here the number of $(a, b)$-transverse intersection points is counted with multiplicity. 

\subsection{Closed geodesics on a complete hyperbolic surface}
\begin{figure}[H]
\centering
\begin{overpic}[scale=0.5]{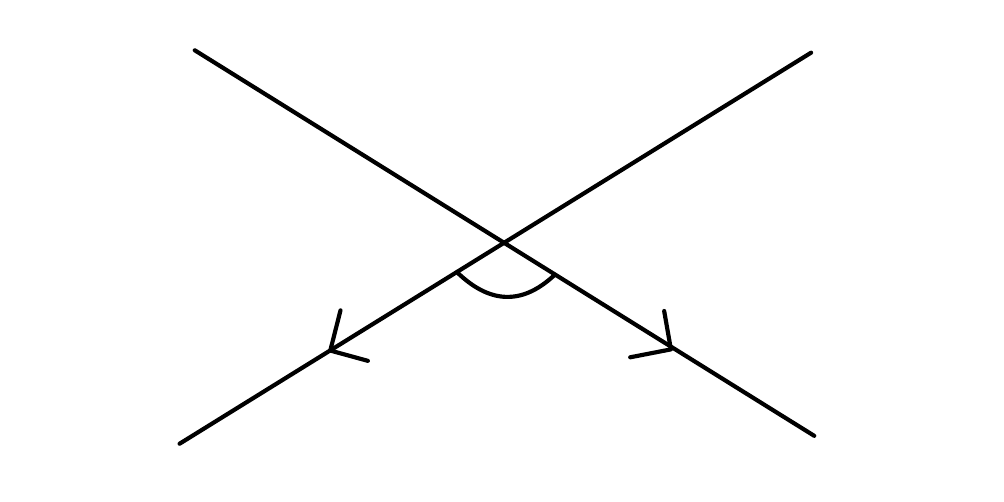}
\put(48,30){$P$}
\put(30,8){$\Al(X)$}
\put(60,8){$\Be(X)$}
\put(45,15){$\phi_{P}(X)$}
\end{overpic}
\caption{The forward angle $\phi_{P}(X)$}
\label{fig.forward_angle_by_metric}
\end{figure}
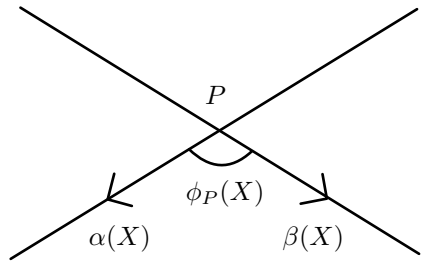

We assume that the Euler characteristic of \( \Sigma \) is negative so that \( \Sigma \) admits a complete hyperbolic metric which makes us identify the universal covering space of \( \Sigma \) with the upper half-plane \( \mathbb{H} \). Denote by $\mathrm{pr} \colon \mathbb{H} \to \Sigma$ the universal covering map. We identify the universal cover of \( \Sigma \) with the upper half-plane $\mathbb{H}$. Let $\alpha, \beta \in \hat\pi$. We denote by $\Al(X)$ the geodesic representative of $\Al$ with respect to a complete hyperbolic metric $X$. For each \( (\Al(X) ,\Be(X) )\)-transverse intersection point $P$, \emph{the X-forward angle of $\Al$ and $\Be$ at $P$}, denoted by $\phi_{P}(X)$, is defined as the angle between the directions of \( \Al(X) \) and \( \Be(X) \) at \( P \), where \( 0 < \phi_{P}(X) < \pi \). We also define the forward angle $\phi_P(X)$ at a self-intersection point $P$ of \( \Al(X) \) in the same way. In addition, for two directed geodesics in the upper half-plane $\mathbb{H}$ and their intersection point $P$, we also define the forward angle $\phi_P$ in the same way. See Figure \ref{fig.forward_angle_by_metric}.

For each $\Al \in \hat{\pi}$, denote by $\ell_{\Al}(X)$ the length of $\Al(X)$. To simplify notation, we write $\ell_{\Al}$ instead of $\ell_{\Al}(X)$ if there is no confusion and similarly for $\phi_{P}(X)$.

\subsection{Hyperbolic isometries}

Since any deck transformation of the universal covering map $\mathrm{pr} \colon \mathbb{H} \to \Sigma$ is an isometry of $\mathbb{H}$, we review some basic properties of hyperbolic isometries. In the hyperbolic plane \(\mathbb{H}\), we write the hyperbolic distance between two points as \( d(\cdot, \cdot) \) throughout this paper.  
For a hyperbolic isometry \( g \), the axis \( A_g \) is the geodesic connecting the two fixed points of \( g \), and the translation length \( t_g \) is the distance that \( g \) moves each point on \( A_g \). The following propositions are direct consequences of the properties of hyperbolic isometries.

\begin{prop} \label{hyperbolic_isometry}
Let $g$ be an isometry of the hyperbolic plane $\mathbb{H}$. Then the following are equivalent,
\begin{enumerate}[\upshape(i)]
\item $g$ is hyperbolic.
\item $g = \rho_2 \rho_1$ where $\rho_1,\rho_2$ are reflections about some disjoint geodesics $L_1,L_2$.
\item $g = r_2 r_1$ where $r_1,r_2$ are rotations of order two about some distinct points $v_1,v_2$.
\end{enumerate}
Then the axis is orthogonal to $L_1,L_2$, and passes through $v_1,v_2$. Its direction is from $L_1$ to $L_2$ and from $v_1$ to $v_2$. Moreover $d(L_1,L_2)=d(v_1,v_2)=\frac{t_g}{2}$
\end{prop}

For each transverse intersection point $P \in \alpha(X) \cap \beta(X)$, we will study how the lifts of $\alpha(X)$, $\beta(X)$, and $|\alpha_P \beta_P|(X)$ are positioned in the upper half-plane $\mathbb{H}$ in the next section. For this purpose, we briefly recall how the following theorem can be derived by using Proposition \ref{hyperbolic_isometry}.

\begin{lemma} \label{thm_cosh}
\cite[Theorem 7.38.6]{Beardon1983}
Let $g$ and $h$ be hyperbolic isometries of the hyperbolic plane and suppose that the axes $A_{g}$ and $A_{h}$ intersect at a point $P$. Denote by $\theta$ the forward angle of $A_{g}$ and $A_{h}$ at $P$. Then the product $gh$ is hyperbolic and 
\begin{align}
\cosh\left(\frac{t_{gh}}{2}\right) = \cosh\left(\frac{t_{g}}{2}\right)\cosh\left(\frac{t_{h}}{2}\right) + \sinh\left(\frac{t_{g}}{2}\right)\sinh\left(\frac{t_{h}}{2}\right)\cos\theta.
\end{align}

\begin{figure}[H]
\centering
\begin{overpic}[scale=0.58]{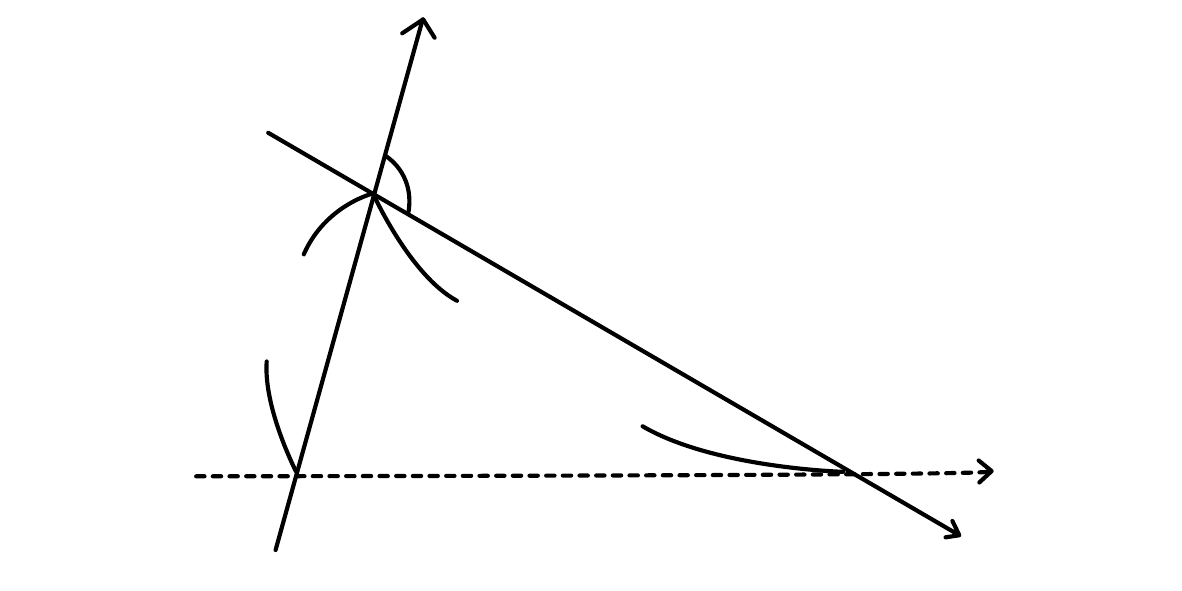}
\put(29,36){$P$}
\put(26,7){$v_1$}
\put(70,7){$v_2$}
\put(23,24){$\frac{t_h}{2}$}
\put(45,18){$\frac{t_g}{2}$}
\put(82,3){$A_{g}$}
\put(38,47){$A_{h}$}
\put(85,11){$A_{gh}$}
\put(35.5,35){$\theta$}
\end{overpic}
\caption{Relative position of the axes \( A_g \), \( A_h \), and \( A_{gh} \) for hyperbolic isometries \( g \), \( h \), and \( gh \).}
\label{fig_thm_cos}
\end{figure}
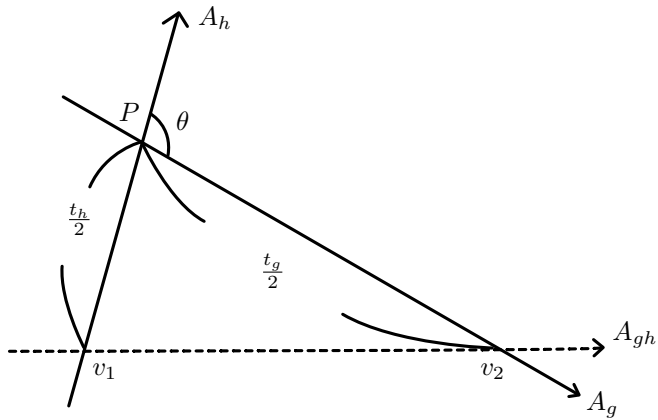
\end{lemma}

\begin{proof}
Let \( v_1 \) be the point on the axis \( A_h \) located at a distance \( \frac{t_h}{2} \) backward along its direction, and \( v_2 \) the point on the axis \( A_g \) located at a distance \( \frac{t_g}{2} \) forward along its direction (see Figure~\ref{fig_thm_cos}).

By Proposition~\ref{hyperbolic_isometry}, we have $g = r_2 r_P$ and $h = r_P r_1$ where \( r_1 \) and \( r_2 \) are rotations of order two about the points \( v_1 \) and \( v_2 \), respectively, and \( r_P \) is a rotation of order two about the intersection point \( P \). Therefore, we obtain
\[
gh = (r_2 r_P)(r_P r_1) = r_2 r_1.
\]
By Proposition~\ref{hyperbolic_isometry}, the product \( gh = r_2 r_1 \) is a hyperbolic isometry whose axis \( A_{gh} \) is the geodesic \( L \) joining the points \( v_1 \) and \( v_2 \), whose direction is from \( v_1 \) to \( v_2 \), and whose translation length satisfies \( d(v_1, v_2) = \frac{t_{gh}}{2} \).
By the hyperbolic law of cosines, we obtain
\[
\cosh\left(\frac{t_{gh}}{2}\right) = \cosh\left(\frac{t_{g}}{2}\right)\cosh\left(\frac{t_{h}}{2}\right) + \sinh\left(\frac{t_{g}}{2}\right)\sinh\left(\frac{t_{h}}{2}\right)\cos \theta.
\]
\end{proof}

\section{Separability criteria}
\subsection{A zigzag curve induced by an intersection point of two closed geodesics}

This section is devoted to the proof of the separability criteria, Theorems \ref{WSC1} and \ref{SSC1} in Introduction.
In order to prove them, we analyze a zigzag curve induced by an intersection of two geodesics.
Let \(\alpha, \beta \in \hat{\pi}\), and \(P\) an \((\alpha(X), \beta(X))\)-transverse intersection point with respect to a complete hyperbolic metric \(X\) on \(\Sigma\). We denote by $\Al(X)_P$ and $\Be(X)_P$ the geodesics $\Al(X)$ and $\Be(X)$ regarded as starting at $P$. Each term $|\alpha_P\beta_P|$ in the Goldman bracket is represented by the loop obtained by first traversing $\Al(X)_P$, followed by $\Be(X)_P$. Kabiraj [\cite{Kabiraj2016}, Figure 1, p.2844] studied lifts of the geodesic representative $|\alpha_P\beta_P|(X)$ in the upper half-plane $\mathbb{H}$ in the case $\alpha$ is simple, to compute the center of the Goldman Lie algebra of surfaces. In computing the center of the Thurston--Wolpert--Goldman Lie algebra, a Lie subalgebra of the Goldman Lie algebra, Chas and Kabiraj \cite{Chas-Kabiraj2022} studied such lifts in detail and called the resulting bi-infinite piecewise geodesic a \emph{zigzag curve}. Our proof also relies on a detailed analysis of zigzag curves, including the case where $\alpha$ is not necessarily simple.

We fix a lift \( P_0' \) of the intersection point \( P \) to the universal cover \( \mathbb{H} \). Let \( P_0'' \) be the endpoint of the lift of \( \Al(X)_P \) that starts at \( P_0' \). Then let \( P_1' \) be the endpoint of the lift of \( \Be(X)_P \) that starts at \( P_0'' \), and let \( P_1'' \) be the endpoint of the lift of \( \Al(X)_P \) that starts at \( P_1' \). Continuing this process inductively, we define a sequence of points \( P_i' \) and \( P_i'' \) for all \( i \geq 0 \). In the same way, we define \( P_i' \) and \( P_i'' \) for all \( i < 0 \) by proceeding in the opposite direction. By connecting these alternating lifts of \( \Al(X)_P \) and \( \Be(X)_P \), we obtain a bi-infinite piecewise geodesic, and we call it a \emph{zigzag curve} \( C \) induced by the \((\alpha(X), \beta(X))\)-transverse intersection point \( P \). For each integer \( i \), let \( M_{2i} \) be the midpoint of the segment \( P_i'P_i'' \), and let \( M_{2i+1} \) be the midpoint of the segment \( P_i''P_{i+1}' \). By Theorem~\ref{thm_cosh}, all the points \( M_i \) lie on a common geodesic. We define \( L \) to be this geodesic, oriented in the same direction as \( C \).

\begin{figure}[H]
\centering
\begin{overpic}[scale=0.51]{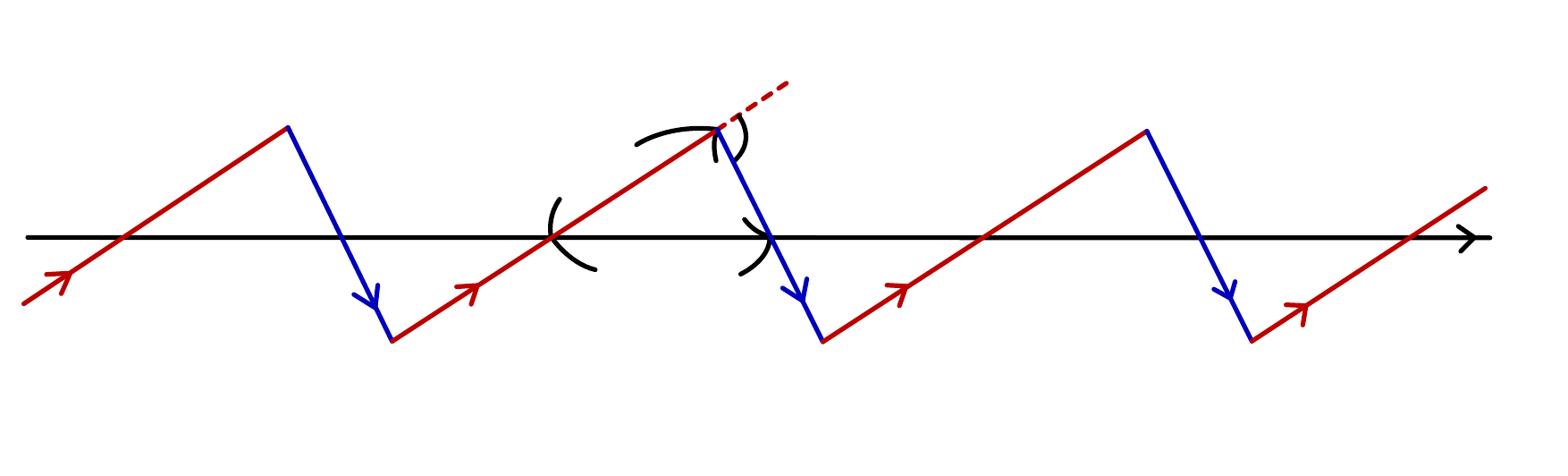}
\put(17,23){$P'_{-1}$}
\put(43,23){$P'_{0}$}
\put(72,23){$P'_{1}$}
\put(24,6){$P''_{-1}$}
\put(52,6){$P''_{0}$}
\put(80,6){$P''_{1}$}
\put(48.5,21){$\phi_{P}$}
\put(38.7,11){$\frac{\ell_{\left| \Al_{P} \Be_{P} \right|}}{2}$}
\put(44,16.9){$\frac{\ell_\Al}{2}$}
\put(37,20){$\frac{\ell_\Be}{2}$}
\put(22,16){$M_{-2}$}
\put(50,16){$M_{0}$}
\put(77,16){$M_{2}$}
\put(7,12){$M_{-3}$}
\put(33,12){$M_{-1}$}
\put(62,12){$M_{1}$}
\put(89,12){$M_{3}$}
\put(96,14){$L$}
\put(96,18){$C$}
\end{overpic}
\caption{A zigzag curve from alternating lifts of \(\Al(X)_P\) (blue) and \(\Be(X)_P\) (red).}
\label{fig_zigzag0}
\end{figure}
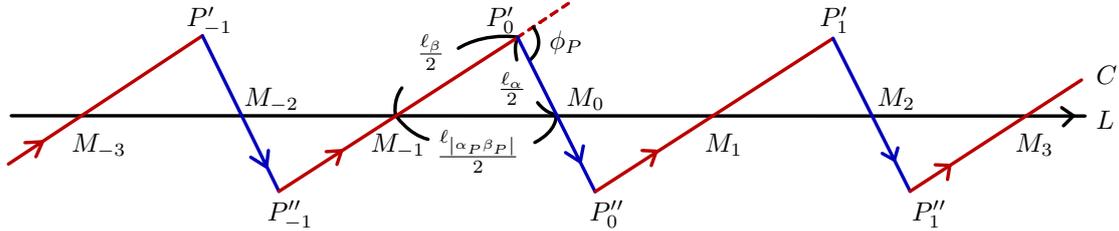

Consider the loop obtained by first traveling along \(\Al(X)_P\) and then along \(\Be(X)_P\). This loop is freely homotopic to \(|\Al_P \Be_P|\). This observation leads to the following lemma.

\begin{lemma}\label{thm_cosh_loop}
Let \( X \) be a complete hyperbolic metric. Let \( \Al, \Be \in \hat{\pi} \) with \( P \in \Al(X) \cap \Be(X) \). Then we have
\[
\cosh\left( \frac{\ell_{\left| \Al_{P} \Be_{P} \right|}}{2} \right) = \cosh\left( \frac{\ell_\Al}{2} \right)\cosh\left( \frac{\ell_\Be}{2} \right) + \sinh\left( \frac{\ell_\Al}{2} \right)\sinh\left( \frac{\ell_\Be}{2} \right)\cos\phi_P.
\]
\end{lemma}

\subsection{Intersection of two zigzag curves}\label{IOTZC}

In order to investigate how a zigzag curve \( C \) and its reflection with respect to a geodesic orthogonal to \( L \) are configured, we now introduce the following notation. We write \( C^{-1} \) for the zigzag curve obtained from \( C \) by reversing the direction of each segment. Let \( \theta_0 \) denote the forward angle of \( L \) and \( P_0'P_0'' \) at \( M_0 \). For each integer \( i \in \mathbb{Z} \), let \( H_i \) be the foot of the perpendicular from \( P_i' \) to \( L \) and \( a \in \mathbb{R}_{\geq0} \) the length of the segment \( M_0 H_0 \). For $u \in \mathbb{R}$, let $\mathcal{U}_u$ be the geodesic perpendicular to $L$ passing through the point on the directed geodesic $L$ at signed distance $u$ from $H_0$. Let \( \rho_{\mathcal{U}_u} \) denote the reflection with respect to \( \mathcal{U}_u \). We consider \( \rho_{\mathcal{U}_u}(C) \), obtained by reflecting each segment of \( C \). Each segment of \( \rho_{\mathcal{U}_u}(C) \) carries the direction induced from \( C \) by the reflection. Define the zigzag curve \( D_u = \rho_{\mathcal{U}_u}(C)^{-1} \). For each integer \( i \in \mathbb{Z} \), let \( Q_i' = \rho_{\mathcal{U}_u}(P_{-i}'') \), \( Q_i'' = \rho_{\mathcal{U}_u}(P_{-i}') \), and \( N_i = \rho_{\mathcal{U}_u}(M_{-i}) \). Define \( \mathcal{V}_u \) to be the geodesic perpendicular to \( L \) obtained by translating \( \mathcal{U}_u \) by \( \frac{1}{2} \ell_{\left| \Al_{P} \Be_{P} \right|} \) in the direction of \( L \). Let \( \rho_{\mathcal{V}_u} \) denote the reflection with respect to \( \mathcal{V}_u \).

\begin{figure}[H]
\centering
\begin{overpic}[scale=0.68]{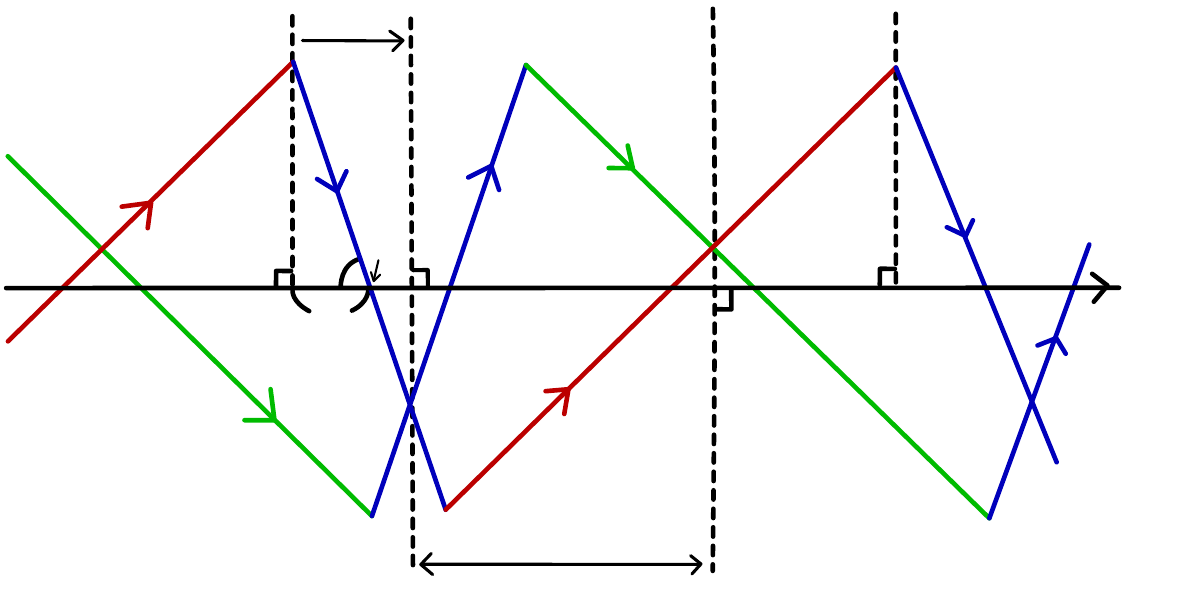}
\put(21,45){$P'_{0}$}
\put(72.5,45){$P'_{1}$}
\put(30,3.5){$Q'_{0}$}
\put(81,3.5){$Q'_{1}$}
\put(36.7,3.5){$P''_{0}$}
\put(42,46){$Q''_{0}$}
\put(12,26.5){$N_{-1}$}
\put(5,23){$M_{-1}$}
\put(38,23){$N_{0}$}
\put(30.8,28.8){$M_{0}$}
\put(64,26.5){$N_{1}$}
\put(56,23){$M_{1}$}
\put(87.7,27){$N_{2}$}
\put(83.5,27){$M_{2}$}
\put(29,47){$u$}
\put(27,28.5){$\theta_0$}
\put(27.5,22.5){$a$}
\put(22,22.5){$H_0$}
\put(73,23){$H_1$}
\put(34,49.5){$\mathcal{U}_u$}
\put(59.5,50){$\mathcal{V}_u$}
\put(44,4){$\frac{\ell_{\left| \Al_{P} \Be_{P} \right|}}{2}$}
\put(89,8){$C$}
\put(92,30){$D_u$}
\put(96,25){$L$}
\end{overpic}
\caption{A zigzag curve \(C\) and its reflection \(D_u\) across \(\mathcal{U}_u\), together with the geodesics \(L\), \(\mathcal{U}_u\), \(\mathcal{V}_u\) and the points \(M_i, N_i, P_i', P_i'', Q_i', Q_i'', H_i\).}
\label{fig_zigzag_def2}
\end{figure}
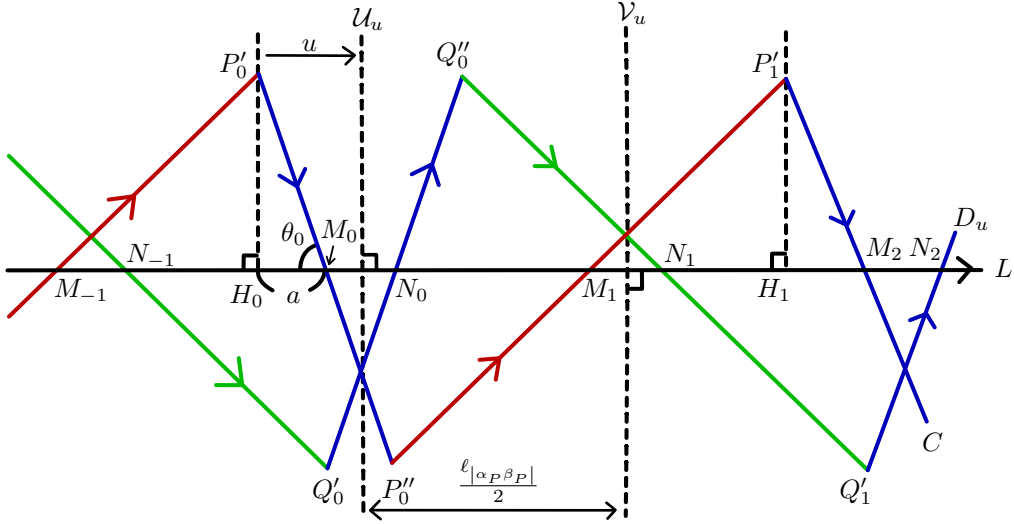

By Proposition \ref{hyperbolic_isometry}, the composition \( \rho_{\mathcal{V}_u} \circ \rho_{\mathcal{U}_u} \) is the hyperbolic isometry whose axis is \( L \), whose direction agrees with the direction of \( L \), and whose translation length is \( \ell_{\left| \Al_{P} \Be_{P} \right|} \). Hence, we get 
\[
\rho_{\mathcal{V}_u}(C) = (\rho_{\mathcal{V}_u} \circ \rho_{\mathcal{U}_u}) \circ \rho_{\mathcal{U}_u}(C) = (\rho_{\mathcal{V}_u} \circ \rho_{\mathcal{U}_u})(D_u^{-1}) = D_u^{-1},
\]
This follows from the identity \( \rho_{\mathcal{U}_u}^2 = \mathrm{id} \), together with the observation that \( \rho_{\mathcal{V}_u} \circ \rho_{\mathcal{U}_u} \) maps each \( M_i \) to \( M_{i+2} \), and hence each \( N_i \) to \( N_{i+2} \).

To classify the relative position of $C$ and $D_u$, we vary $u$ while keeping $C$ fixed, which makes $D_u$ move along $L$. Since translating $D_u$ by an integer multiple of $\ell_{|\Al_P \Be_P|}$ yields the same configuration, it suffices to consider $D_u$ within a single period of length $\ell_{|\Al_P \Be_P|}$. In addition, for \(u, u' \in \mathbb{R}\), \(D_{u'}\) is obtained from \(D_u\) by translating it along \(L\) by distance \(2(u' - u)\), since \(D_u\) (resp. \(D_{u'}\)) is obtained from \(C\) by the reflection across \(\mathcal{U}_u\) (resp. \(\mathcal{U}_{u'}\)). Hence, $D_{u+\frac{1}{2}\,\ell_{|\Al_P \Be_P|}}$ is obtained from $D_u$ by translating it along $L$ by distance $\ell_{|\Al_P \Be_P|}$. Therefore, we may consider only the case $0<u \leq \frac{1}{2}\,\ell_{|\Al_P \Be_P|}$ in the following. Here we don't take the interval $0\leq u < \frac{1}{2}\,\ell_{|\Al_P \Be_P|}$ to avoid index shifts in the following lemma.

\begin{lemma}\label{zigzag_lem1}
Suppose that $\ell_{\Al} \leq \ell_{\Be}$ and \( 0 < u \leq \frac{1}{2} \ell_{|\Al_P \Be_P|} \). Then the segments \( Q_0''Q_1' \) and \( P_0''P_1' \) intersect transversely.
\end{lemma}

\begin{proof}
It suffices to show that the geodesic \( \mathcal{V}_u \) intersects the segment \( P_0''P_1' \) since $\rho_{\mathcal{V}_u}(Q_0''Q_1')=P_1'P_0''$ (See Figure \ref{fig_zigzag_def2}). We consider two cases depending on whether \( \theta_0 < \frac{\pi}{2} \) or \( \theta_0 \geq \frac{\pi}{2} \).

In the case \( \theta_0 < \frac{\pi}{2} \), the assumption \( \ell_\Al \leq \ell_\Be \) implies that \(\angle P_0' M_{-1} M_{0} \leq \theta_0 < \frac{\pi}{2}\) by the hyperbolic law of sines. Since both angles \( \angle P_0' M_{-1} M_{0} \) and \( \theta_0 \) are less than \( \frac{\pi}{2} \), the sum of the interior angles of either \( \triangle P_0' H_0 M_{0} \) or \( \triangle P_0' H_0 M_{-1} \) exceeds \( \pi \) if \( H_0 \) lies outside the segment \( M_{0}M_{-1} \). This contradicts the Gauss--Bonnet theorem, so \( H_0 \) lies on the segment \( M_{0}M_{-1} \). Again, by the same assumption \(\ell_\Al \leq \ell_\Be\), the length of the segment \(M_{0}H_{0}\) is smaller than that of \(M_{-1}H_{0}\) by the hyperbolic law of cosines. Hence, we obtain $2a \leq \tfrac{1}{2}\ell_{|\Al_P \Be_P|}$ (See Figure \ref{SC15}). In addition, the condition \(0 < u \leq \tfrac{1}{2}\ell_{|\Al_P \Be_P|}\) gives $\tfrac{1}{2}\ell_{|\Al_P \Be_P|} < u + \tfrac{1}{2}\ell_{|\Al_P \Be_P|} \leq \ell_{|\Al_P \Be_P|}.$ Combining these inequalities, we have
\[
2a < u + \frac{1}{2}\ell_{|\Al_P \Be_P|} \leq \ell_{|\Al_P \Be_P|},
\]
which means that, in terms of the signed distance from $H_0$ along $L$, the intersection point of the geodesic $\mathcal{V}_u$ with $L$ lies between the foot of the perpendicular from $P_0''$ to $L$ and the point $H_1$. Therefore, $\mathcal{V}_u$ intersects the segment $P_0''P_1'$.

\begin{figure}[H]
\centering
\begin{overpic}[scale=0.85]{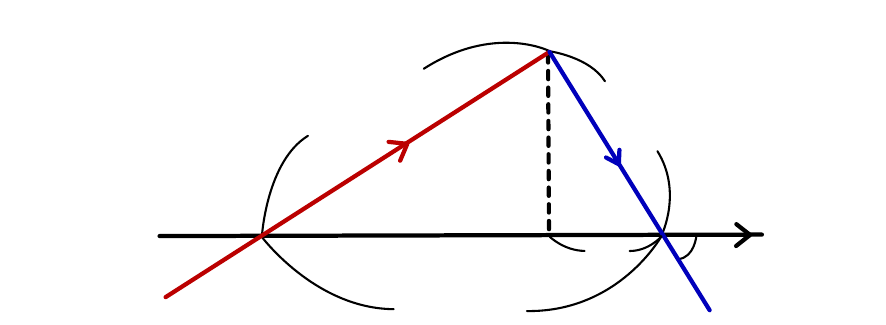}
\put(61,32){$P_0'$}
\put(59,7){$H_0$}
\put(75.5,11.3){$M_0$}
\put(22,11.3){$M_{-1}$}
\put(67.5,7){$a$}
\put(47,1){$\frac{1}{2}\ell_{|\Al_P \Be_P|}$}
\put(40,26){$\frac{1}{2}\ell_{\Be}$}
\put(71,24){$\frac{1}{2}\ell_{\Al}$}
\put(79,6){$\theta_0$}
\put(81,0){$C$}
\put(87,9){$L$}
\end{overpic}
\caption{Proof of Lemma \ref{zigzag_lem1} in the case \( \theta_0 < \frac{\pi}{2} \)}
\label{SC15}
\end{figure}
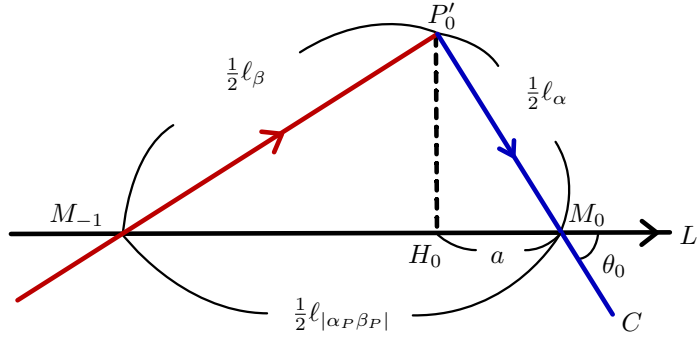

Next, we consider the case \( \theta_0 \geq \frac{\pi}{2} \). By the conditions $0 < u \leq \frac{1}{2}\ell_{|\Al_P \Be_P|} $ and $a\geq0$, we obtain
\[
\frac{1}{2}\ell_{|\Al_P \Be_P|} -a < u + \frac{1}{2}\ell_{|\Al_P \Be_P|} \leq \ell_{|\Al_P \Be_P|},
\]
which means that, in terms of the signed distance from $H_0$ along $L$, the intersection point of the geodesic $\mathcal{V}_u$ with $L$ lies on the segment \( M_1H_1 \) (See Figure \ref{SC14}). Hence, the geodesic \( \mathcal{V}_u \) enters the triangle \( \triangle M_1H_1P_1' \) from the segment \( M_1H_1 \) and must exit through either the segment \( M_1P_1' \) or \( P_1'H_1 \). Since \( \angle M_1H_1P_1' = \frac{\pi}{2} \), the geodesic \( \mathcal{V}_u \) cannot exit through the side \( P_1'H_1 \). Therefore, it exits through the side \( M_1P_1' \), and hence intersects the segment \( P_0''P_1' \). We denote by \( T' \) the transverse intersection point of \( P_0''P_1' \) and \( Q_0''Q_1' \).

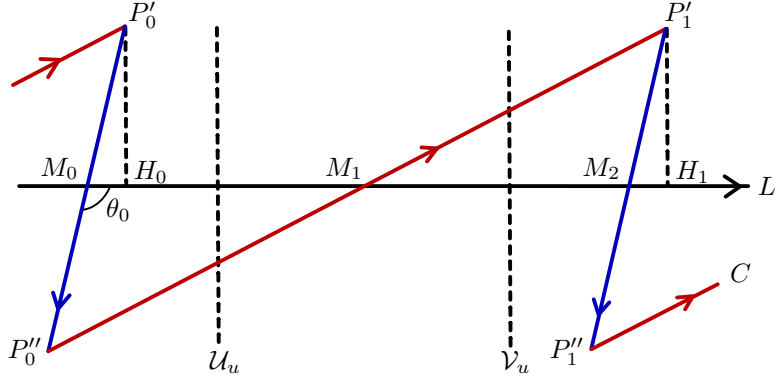
\begin{figure}[H]
\centering
\begin{overpic}[scale=0.72]{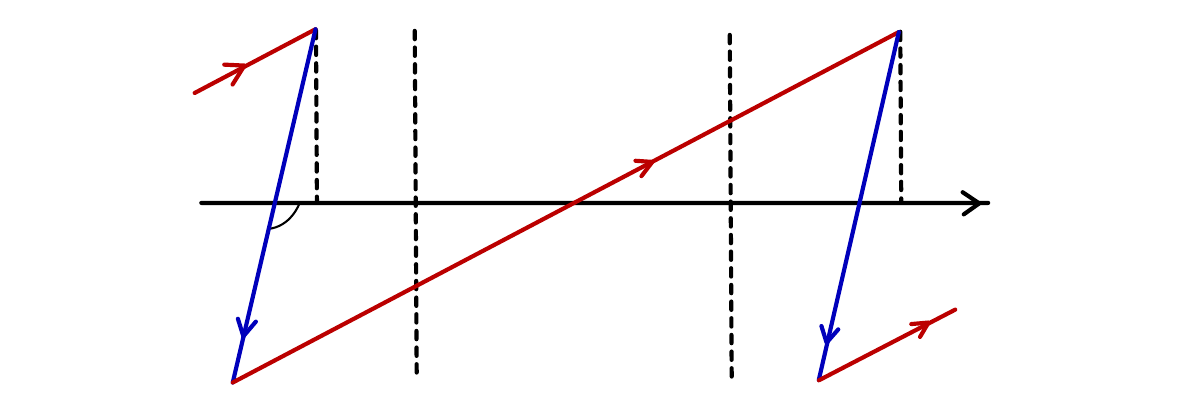}
\put(19,18.8){$M_0$}
\put(45,18.8){$M_1$}
\put(68.5,18.8){$M_2$}
\put(27,33){$P_0'$}
\put(34.5,1){$\mathcal{U}_u$}
\put(61,1){$\mathcal{V}_u$}
\put(27.5,18.5){$H_0$}
\put(16,2){$P_0''$}
\put(76,33){$P_1'$}
\put(77,18.5){$H_1$}
\put(65.5,2){$P_1''$}
\put(84.5,17){$L$}
\put(82,9){$C$}
\put(25,15){$\theta_0$}
\end{overpic}
\caption{Proof of Lemma \ref{zigzag_lem1} in the case \( \theta_0 \geq \frac{\pi}{2} \)}
\label{SC14}
\end{figure}

\end{proof}

In the next lemma, we have to consider several distinct geodesics in \(\mathbb{H}\) intersecting at a single point in some degenerate cases. We recall the notion of an \((a,b)\)\textit{-transverse intersection point} introduced in Subsection~\ref{A generalized notion of intersection points}. In the following, a \((PQ,RS)\)\textit{-transverse intersection point}  is the intersection of the geodesics in \(\mathbb{H}\) containing the segments \(PQ\) and \(RS\).

\begin{lemma}\label{zigzag_lem2}  
Suppose that \( \ell_{\Al} < \ell_{\Be} \), \( 2a \leq u \leq \frac{1}{2} \ell_{|\Al_P \Be_P|} \), and \( \theta_0 < \frac{\pi}{2} \). Then there exist a \((P_0''P_1',\, Q_0'Q_0'')\)-transverse intersection point \( R'\) and a \((P_0'P_0'',\, Q_{-1}''Q_0')\)-transverse intersection point \(\bar{R}'\) such that \( \phi_{R'} = \phi_{\bar{R}'} < \phi_P \).
\end{lemma}

\begin{proof}
We classify the relative positions of the zigzag curves \( C \) and \( D_u \) into four cases: \(\mathrm{(I)}\) \( 0 < u < 2a \), \(\mathrm{(II)}\) \( u = 2a \), \(\mathrm{(III)}\) \( 2a < u < \frac{1}{2} \ell_{\left| \Al_{P} \Be_{P} \right|} \), and \(\mathrm{(IV)}\) \( u = \frac{1}{2} \ell_{\left| \Al_{P} \Be_{P} \right|} \) (See Figure \ref{fig_zigzag5_3}).

\begin{figure}[H]
\centering
\begin{overpic}[width=1.02\linewidth,trim=1cm 0 0 0,clip]{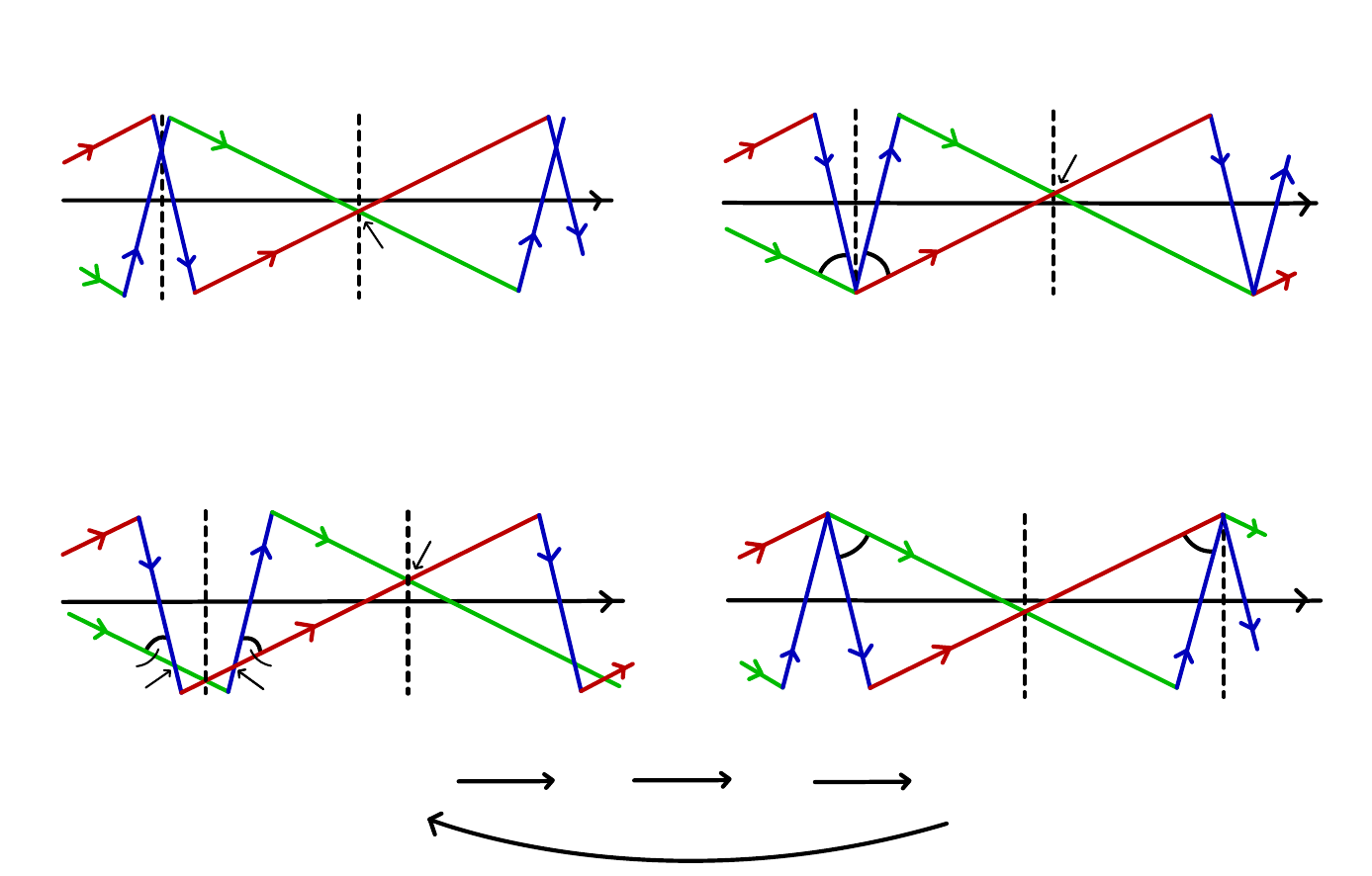}
\put(0,63){$\mathrm{(I)}$ $ 0<u<2a$}
\put(50,63){$\mathrm{(II)}$ $ u=2a$}
\put(0,34){$\mathrm{(III)}$ $ 2a<u<\frac{1}{2}\ell_{|\Al_P\Be_P|}$}
\put(50,34){$\mathrm{(IV)}$ $ u=\frac{1}{2}\ell_{|\Al_P\Be_P|}$}
\put(7.2,42){$\mathcal{U}_u$}
\put(22.2,42){$\mathcal{V}_u$}
\put(39,45){$C$}
\put(39,57.5){$D_u$}
\put(4.5,59){$P_0'$}
\put(35.5,59){$P_1'$}
\put(10,42){$P_0''$}
\put(25,45.5){$T'$}
\put(3.5,42){$Q_0'$}
\put(34.5,42){$Q_1'$}
\put(8.5,59){$Q_0''$}
\put(43,51){$L$}
\put(4,52.4){$N_0$} \put(9,52.4){$M_0$} \put(34,52.4){$N_2$} \put(39,52.4){$M_2$} \put(19.6,52.9){$N_1$} \put(23.3,52.9){$M_1$}

\put(96.5,50.5){$L$}
\put(60,59){$\mathcal{U}_u$}
\put(75,59){$\mathcal{V}_u$}
\put(50.5,42.5){$P_0''=Q_0'=R'=\bar{R}'$}
\put(56,59){$P_0'$}
\put(63,59){$Q_0''$}
\put(87,59){$P_1'$}
\put(86.5,42){$P_1''=Q_1'$}
\put(95,46){$C$}
\put(94,56){$D_u$}
\put(62.5,48){$\phi_{R'}$}
\put(56,48){$\phi_{\bar{R'}}$}
\put(77.5,55.5){$T'$}
\put(63,52){$N_0$} \put(55.3,52){$M_0$} \put(92.5,48.4){$N_2$} \put(85.5,52){$M_2$} \put(76.5,48.6){$N_1$} \put(72.5,48.6){$M_1$}

\put(7.5,11){$P_0''$}
\put(12,11){$Q_0'$}
\put(16.5,15){$\phi_{R'}$}
\put(16,12){$R'$}
\put(4,12){$\bar{R'}$}
\put(2,15){$\phi_{\bar{R'}}$}
\put(5,28.5){$P_0'$}
\put(35,28.5){$P_1'$}
\put(15,28.5){$Q_0''$}
\put(10.5,28.5){$\mathcal{U}_u$}
\put(26,28.5){$\mathcal{V}_u$}
\put(28,26){$T'$}
\put(44,16){$C$}
\put(43,13){$D_u$}
\put(43.5,20){$L$}
\put(15.2,21.5){$N_0$} \put(4,21.5){$M_0$} \put(34.4,21.6){$M_2$} \put(28.7,18.5){$N_1$} \put(22,18.5){$M_1$}

\put(52,28.5){$P_0'=Q_{-1}''=\bar{R}'$}
\put(80,28.5){$P_1'=Q_0''=R'=T'$}
\put(96.5,20){$L$}
\put(61,12){$P_0''$}
\put(54,12){$Q_{-1}'$}
\put(73,11){$\mathcal{U}_u$}
\put(88,11){$\mathcal{V}_u$}
\put(91,15){$C$}
\put(92,25){$D_u$}
\put(83.5,22.7){$\phi_{R'}$}
\put(60.6,22.7){$\phi_{\bar{R'}}$}
\put(53.8,21.7){$N_0$} \put(61,18.8){$M_0$} \put(83.3,18.8){$N_2$} \put(90.6,21.7){$M_2$} \put(70.3,22.3){$N_1$} \put(73.8,22.3){$M_1$}

\put(26.8,6.5){$\mathrm{(I)}$}
\put(39,6.5){$\mathrm{(II)}$}
\put(52.2,6.5){$\mathrm{(III)}$}
\put(66,6.5){$\mathrm{(IV)}$}
\end{overpic}
\caption{Possible configurations of the zigzag curves $C$ and $D_u$ in the case \(\ell_{\Al} < \ell_{\Be}\) and \(\theta_0 < \frac{\pi}{2} \)}
\label{fig_zigzag5_3}
\end{figure}
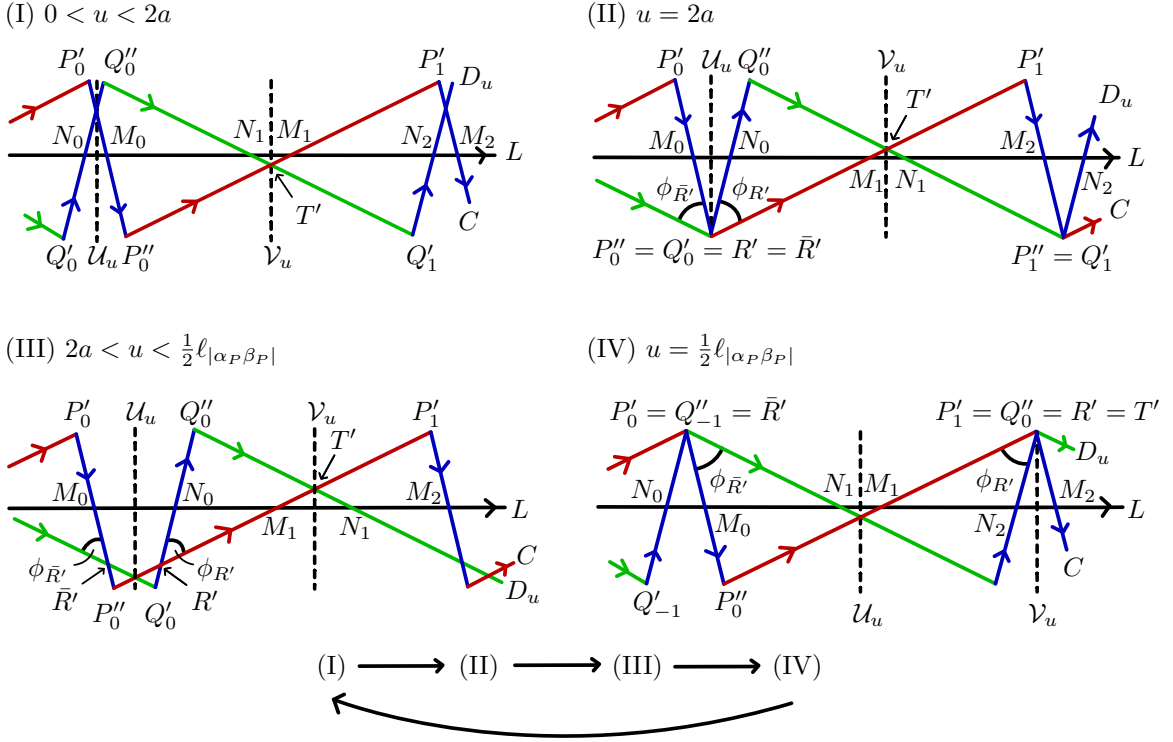

By the assumption \( 2a \leq u \leq \frac{1}{2} \ell_{|\Al_P \Be_P|} \), we only need to consider the cases \(\mathrm{(II)}\), \(\mathrm{(III)}\), and \(\mathrm{(IV)}\).

In Case \(\mathrm{(II)}\), the condition \(u = 2a\) leads to \(P_0'' = Q_0'\), and this point is a \((P_0''P_1', Q_0'Q_0'')\)-transverse intersection point \( R' \) and a \((P_0''P_0', Q_{-1}''Q_0')\)-transverse intersection point \( \bar{R}' \). Moreover, since \( \rho_{\mathcal{U}_u}(C) = D_u^{-1} \), we have \( \phi_R' = \phi_{\bar{R}'} \). Since the sum of the interior angles of the triangle \( T'Q_0''R' \) is less than \( \pi \), we obtain \( \phi_{R'} = \phi_{\bar{R}'} < \phi_P \).

In Case \(\mathrm{(III)}\), by the assumption \(\theta_0 < \frac{\pi}{2}\), we have \(a > 0\). Together with \(2a < u < \frac{1}{2}\,\ell_{|\Al_P \Be_P|}\), we obtain
\[
2a < 2u - 2a < 2u < \ell_{|\Al_P \Be_P|},
\]
which means that, in terms of the signed distance from \( H_0 \) along \( L \), the four distinct perpendiculars from \( P_0'' \), \( Q_0' \), \( Q_0'' \), and \( P_1' \) to \( L \) appear in this order, from closest to farthest from \( H_0 \). In particular, we denote the perpendicular from \( Q_0' \) (resp. \( Q_0'' \)) to \( L \) by \( \mathcal{W}_u' \) (resp. \( \mathcal{W}_u'' \)). Then this inequality implies that the perpendicular \( \mathcal{W}_u' \)(resp. \( \mathcal{W}_u'' \)) intersects the segment \( P_0''P_1' \) transversely at a point \( W' \)(resp. \( W'' \)) (see Figure~\ref{fig_SC17}). Since the point \( N_0 \) is the midpoint of the segment \( Q_0'Q_0'' \), the points \( Q_0' \) and \( Q_0'' \) have the same distance from \( L \) and lie on opposite sides of it. Since both \( W' \) and \( W'' \) lie in the interior of the segment \( P_0''P_1' \), their distances from \( L \) are smaller than those of \( Q_0' \) and \( Q_0'' \). On the closed interval of \( L \) between the two points \( \mathcal{W}_u' \cap L \) and \( \mathcal{W}_u'' \cap L \), both segments \( W'W'' \) and \( Q_0'Q_0'' \) can be regarded as continuous functions of the signed distance from \( L \), where we fix a positive direction perpendicular to \( L \). Moreover, the signed distances of \( Q_0' \) and \( Q_0'' \) from \( L \) have opposite signs, while those of \( W' \) and \( W'' \) have smaller absolute values. Hence, by the intermediate value theorem, the segments \( W'W'' \) and \( Q_0'Q_0'' \) intersect transversely. We denote the intersection point by \( R' \). Similarly, the segments \( P_0'P_0'' \) and \( Q_{-1}'Q_0'' \) intersect transversely at a point \( \bar{R}' \). Again, \( \rho_{\mathcal{U}_u}(C) = D_u^{-1} \) implies \(\phi_{R'} = \phi_{\bar{R}'} < \phi_P \), where the inequality holds for the same reason as in Case \(\mathrm{(II)}\).

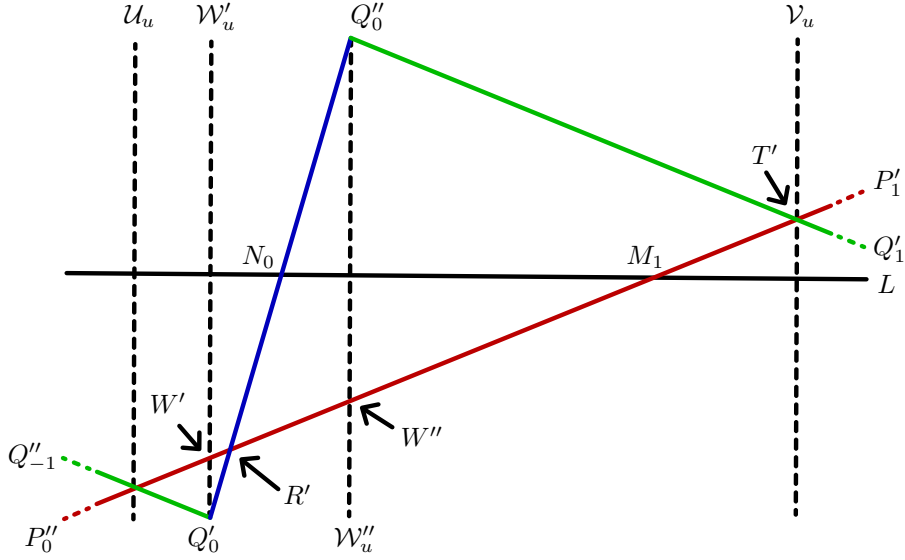
\begin{figure}[H]
\centering
\begin{overpic}[scale=0.55]{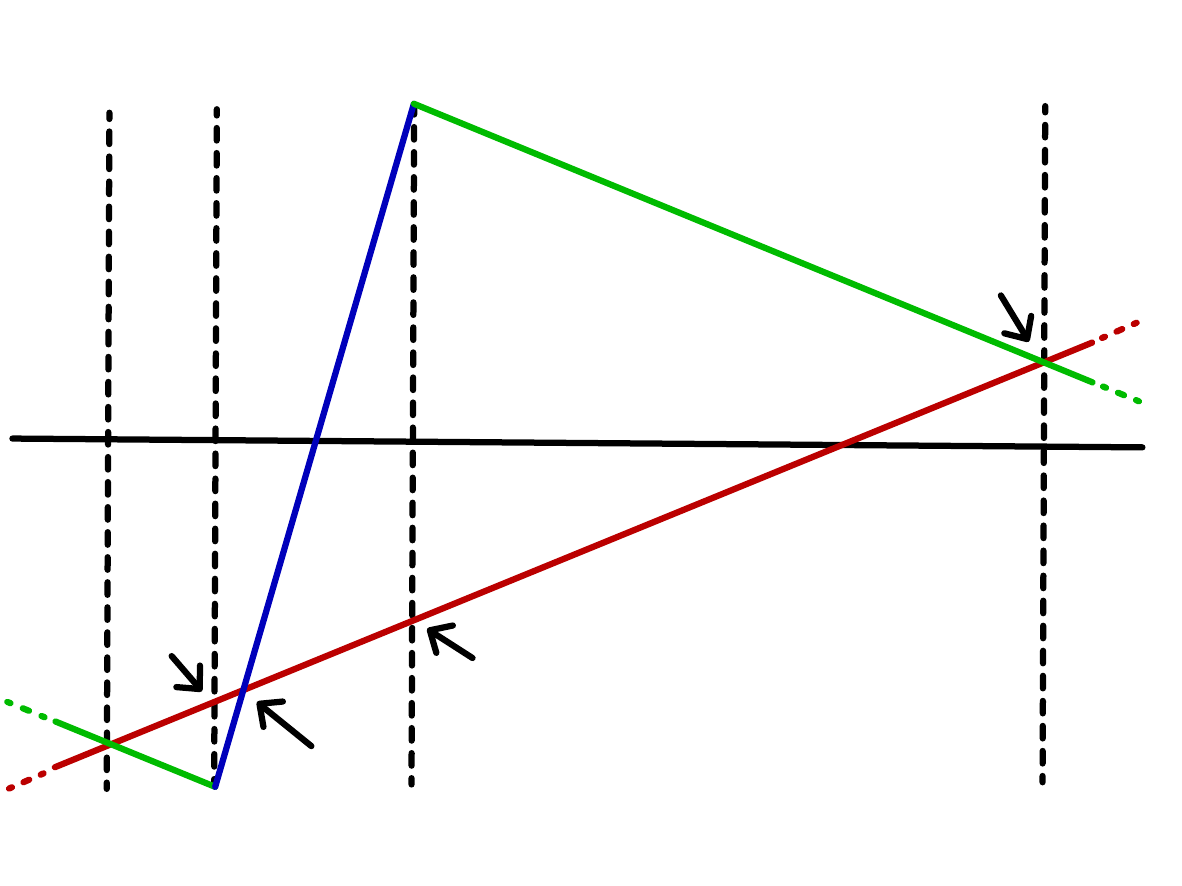}
\put(-10,74){$\mathrm{(III)}$ $ 2a<u<\frac{1}{2}\ell_{|\Al_P\Be_P|}$}
\put(-6,15){$Q_{-1}''$}
\put(-4,5){$P_0''$}
\put(15.5,5){$Q_0'$}
\put(35,68){$Q_0''$}
\put(8,68){$\mathcal{U}_u$}
\put(16.5,68){$\mathcal{W}_u'$}
\put(33,5){$\mathcal{W}_u''$}
\put(87,68){$\mathcal{V}_u$}
\put(83,51){$T'$}
\put(98,35.5){$L$}
\put(97.5,40){$Q_1'$}
\put(97.5,48){$P_1'$}
\put(22,39){$N_0$} \put(68,39){$M_1$}
\put(27,10){$R'$}
\put(11,21){$W'$}
\put(41,17){$W''$}

\end{overpic}
\caption{Case (III) near $\mathcal{U}_u$ and $\mathcal{V}_u$, showing $W'$ and $W''$}
\label{fig_SC17}
\end{figure}

Finally, Case \(\mathrm{(IV)}\) proceeds exactly as in Case \(\mathrm{(II)}\). We omit the details.
\end{proof}

\begin{lemma}\label{zigzag_lem3}  
Suppose that $\ell_{\Al} < \ell_{\Be}$ and \( 0 < u \leq \frac{1}{2} \ell_{|\Al_P \Be_P|} \) and $\theta_0 = \frac{\pi}{2}$. Then there exist a \((P_0'P_1'',\, Q_0'Q_0'')\)-transverse intersection point \( R'\) and a \((P_0'P_0'',\, Q_{-1}''Q_0')\)-transverse intersection point \(\bar{R}'\) such that \( \phi_{R'} = \phi_{\bar{R}'} < \phi_P \).
\end{lemma}

\begin{proof}
We classify the relative positions of the zigzag curves \( C \) and \( D_u \) into two cases: \(\mathrm{(V)}\) \( 0 < u < \frac{1}{2} \ell_{\left| \Al_{P} \Be_{P} \right|}  \), \(\mathrm{(VI)}\) \( u = \frac{1}{2} \ell_{\left| \Al_{P} \Be_{P} \right|}\) (See Figure \ref{fig_zigzag6}).

\begin{figure}[H]
\centering
\begin{overpic}[scale=0.75]{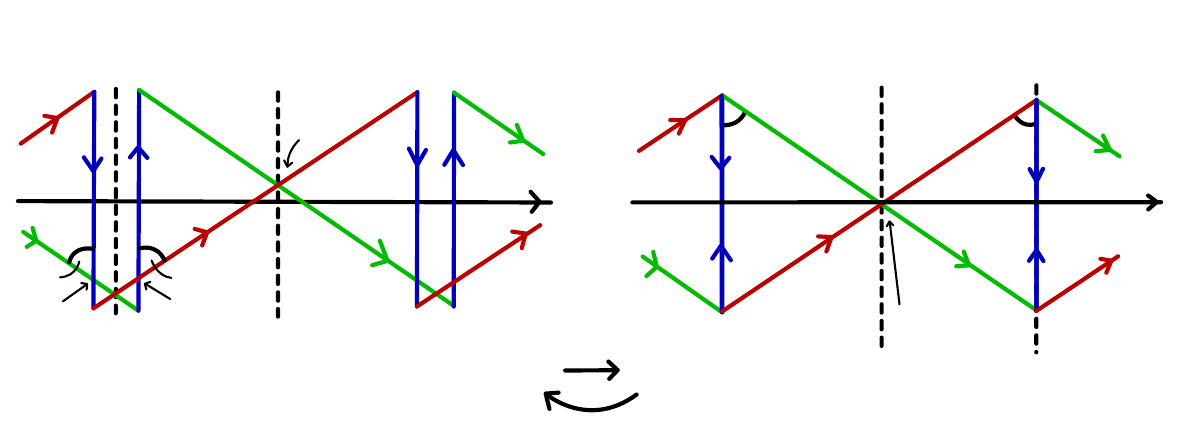}
\put(1,35){$\mathrm{(V)}$ $0 < u < \frac{1}{2} \ell_{\left| \Al_{P} \Be_{P} \right|} $}
\put(52,35){$\mathrm{(VI)}$ $u = \frac{1}{2} \ell_{\left| \Al_{P} \Be_{P} \right|} $}

\put(43.2,5.5){$\mathrm{(V)}$}
\put(53.2,5.5){$\mathrm{(VI)}$}

\put(6,9){$P_0''$}
\put(11,9){$Q_0'$}
\put(15,13.5){$\phi_{R'}$}
\put(14.5,10.5){$R'$}
\put(2.5,10.5){$\bar{R'}$}
\put(1.5,13.5){$\phi_{\bar{R'}}$}
\put(5.5,30.7){$P_0'$}
\put(11.5,30.7){$Q_0''$}
\put(33.5,30.7){$P_1'$}
\put(38,30.7){$Q_1''$}
\put(33.5,9){$P_1''$}
\put(38,9){$Q_1'$}
\put(8.8,31){$\mathcal{U}_u$}
\put(23,31){$\mathcal{V}_u$}
\put(46,17){$C$}
\put(46,23){$D_u$}
\put(47,20){$L$}
\put(25,26){$T'$}
\put(4.2,21.2){$M_0$} \put(12.2,21.2){$N_0$} \put(18.5,21.2){$M_1$} \put(25.5,21.2){$N_1$} \put(31.5,21.2){$M_2$} \put(38.7,21.2){$N_2$} 

\put(54.5,30.5){$P_0'=Q_{-1}'''=\bar{R}'$}
\put(57,9){$P_0''=Q_{-1}'$}
\put(78,31){$P_1'=Q_{0}''=T'=R'$}
\put(89,9){$P_1''=Q_{0}'$}
\put(74,31){$\mathcal{U}_u$}
\put(87,5.5){$\mathcal{V}_u$}
\put(62,24.7){$\phi_{\bar{R'}}$}
\put(84,24.7){$\phi_{R'}$}
\put(95.5,15.5){$C$}
\put(95.5,23){$D_u$}
\put(98.8,20){$L$}
\put(61.5,21){$M_0=N_0$} \put(75,9){$M_1=N_1$} \put(88,18){$M_2=N_2$}

\end{overpic}
\caption{Possible configurations of the zigzag curves $C$ and $D_u$ in the case \(\ell_{\Al} < \ell_{\Be}\) and \(\theta_0 = \frac{\pi}{2} \)}
\label{fig_zigzag6}
\end{figure}
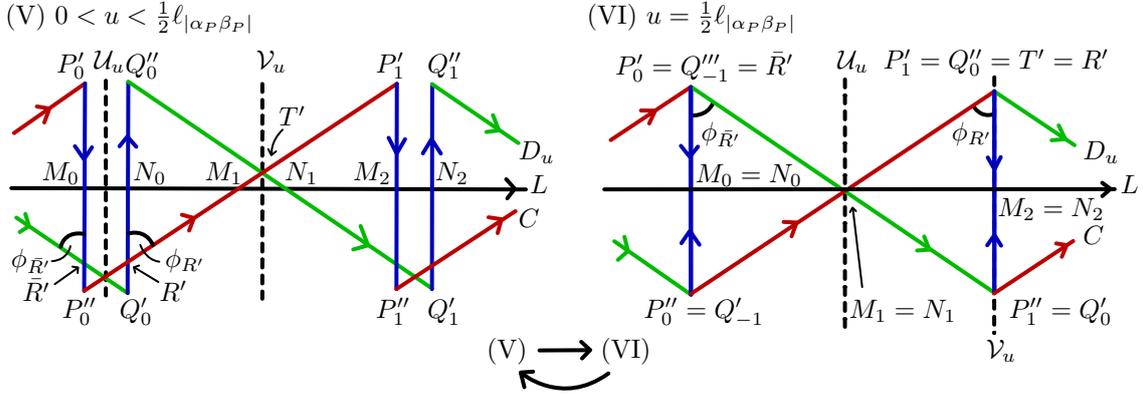

In Case \(\mathrm{(V)}\), since the segment \( Q_0'Q_0'' \) enters either the triangle \( \triangle M_0P_0''M_1 \) or \( \triangle M_1P_1'M_2 \), it must exit from one of these triangles. By the Gauss--Bonnet theorem, it must exit through the segment \( P_1'P_0'' \). Let this intersection point be \( R' \). Similarly, the segment \( P_0'P_0'' \) intersects \( Q_{-1}''Q_0' \), and we denote the intersection point by \( \bar{R}' \). Moreover, since \( \rho_{\mathcal{U}_u}(C) = D_u^{-1} \), we have \( \phi_{R'} = \phi_{\bar{R}'} \). Since the sum of the interior angles of the triangle \( T'Q_0''R' \) is less than \( \pi \), we obtain \( \phi_{R'} = \phi_{\bar{R}'} < \phi_P \).

In Case \(\mathrm{(VI)}\), the configuration is degenerate, and the claim follows in the same way as Case \(\mathrm{(II)}\) of Lemma \ref{zigzag_lem2}. This completes the proof.

\end{proof}

\begin{lemma}\label{zigzag_lem4}
Suppose that $\ell_{\Al} < \ell_{\Be}$ and \( 0 < u \leq \frac{1}{2} \ell_{|\Al_P \Be_P|} \) and $\theta_0 > \frac{\pi}{2}$. Then there exist a \((P_0'P_1'',\, Q_0'Q_0'')\)-transverse intersection point \( R'\) and a \((P_0'P_0'',\, Q_{-1}''Q_0')\)-transverse intersection point \(\bar{R}'\) such that \( \phi_{R'} = \phi_{\bar{R}'} < \phi_P \).
\end{lemma}

\begin{proof}

We classify the relative position of the zigzag curves \( C \) and \( D_u \) into a single case: \(\mathrm{(VII)}\) (See Figure \ref{fig_zigzag7}).

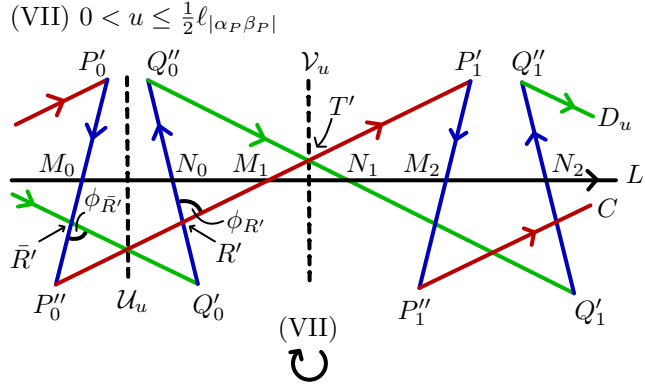
\begin{figure}[H]
\centering
\begin{overpic}[scale=0.75]{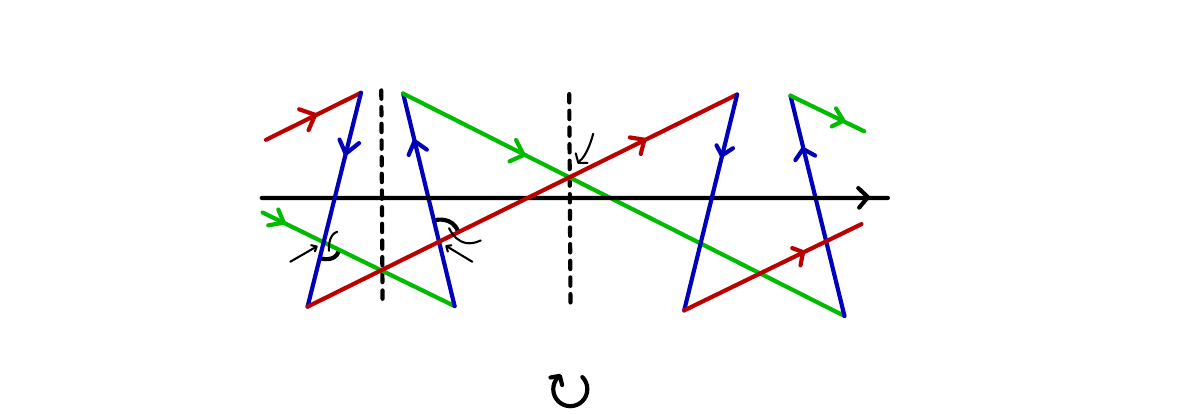}
\put(22,32){\(\mathrm{(VII)}\) \( 0 < u \leq \frac{1}{2} \ell_{|\Al_P \Be_P|}\)}
\put(45.5,4.5){$\mathrm{(VII)}$}

\put(24,6.5){$P_0''$}
\put(38,6.5){$Q_0'$}
\put(41,14){$\phi_{R'}$}
\put(40,11){$R'$}
\put(22,10.5){$\bar{R'}$}
\put(28.5,16){$\phi_{\bar{R'}}$}
\put(28,28){$P_0'$}
\put(34,28){$Q_0''$}
\put(61,28){$P_1'$}
\put(66,28){$Q_1''$}
\put(56,6){$P_1''$}
\put(71.5,6){$Q_1'$}
\put(31.5,7){$\mathcal{U}_u$}
\put(47.5,28){$\mathcal{V}_u$}
\put(73.5,15){$C$}
\put(73.5,23){$D_u$}
\put(76,18){$L$}
\put(50,24){$T'$}
\put(24.5,19){$M_0$} \put(36.5,19){$N_0$} \put(41.5,19){$M_1$} \put(51.5,19){$N_1$} \put(56.5,19){$M_2$} \put(69.5,19){$N_2$} 

\end{overpic}
\caption{The configuration of the zigzag curves $C$ and $D_u$ in the case \(\ell_{\Al} < \ell_{\Be}\) and \(\theta_0 > \frac{\pi}{2} \)}
\label{fig_zigzag7} 
\end{figure}

\noindent The argument proceeds exactly as in Case \(\mathrm{(V)}\) of Lemma \ref{zigzag_lem3}, and we omit the details.

\end{proof}

\subsection{Key results}
For each positive integer \( m \), we continue to denote by \( x^{m} \in \hat{\pi} \) the free homotopy class defined in the Introduction. Note that, at each $(x,y)$-transverse intersection point $P$, $m$ $(x^m,y)$-transverse intersection points are obtained, although on the surface they all correspond to the same point $P$. For each hyperbolic metric $X$, at an $(x(X),y(X))$-transverse intersection point $P$, the forward angle formed by $x(X)$ and $y(X)$ coincides with the $m$ forward angles formed by $x^m(X)$ and $y(X)$ at the corresponding $(x^m(X),y(X))$-transverse intersection points. Hence we denote all these forward angles by $\phi_P(X)$.

One of the key results of this paper is stated in the following lemma. The assumption that \(x^m\) goes around \(x\) at least twice, namely \(m \ge 2\), plays a crucial role in Lemma \ref{xyz_-}. When \(m=1\), that is, when \(x^m\) goes around \(x\) only once, the conclusion does not hold in general.

\begin{lemma} \label{xyz_-}
Let \( X \) be a complete hyperbolic metric, and let \( x, y, z \in \hat{\pi} \) be essential loops with \( \ell_y(X) = \ell_z(X) \). Suppose that \( P \) is an \((x(X),\, y(X))\)-transverse intersection point and \( Q \) is an \((x(X),\, z(X))\)-transverse intersection point satisfying \( \varepsilon_{P}(x,y) = -\varepsilon_{Q}(x,z) \). If there exists a positive integer \( m \geq 2\) such that \( \left| x_P^m y_P \right| = \left| x_Q^m z_Q \right| \). Then there exist an \((x(X),\, y(X))\)-transverse intersection point \( R \) and an \((x(X),\, z(X))\)-transverse intersection point  \( \bar{R} \) such that  
\[
\phi_R = \phi_{\bar{R}} < \phi_P = \phi_Q, 
\quad \text{or} \quad y = z = x^m.
\]
\end{lemma}

\begin{proof}
By Lemma \ref{thm_cosh_loop} and \( \ell_y = \ell_z \), we have $\phi_P = \phi_Q$. Since \(|x^{m}_{P} y_{P}| = |x^{m}_{Q} z_{Q}|\), their common geodesic representative \(|x^{m}_{P} y_{P}|(X)= |x^{m}_{Q} z_{Q}|(X)\) admits a bi-infinite lift to the universal cover \( \mathbb{H} \); denote this lift by \( L \). Along this geodesic line \( L \), there exist two bi-infinite piecewise geodesics \( C \) and \( D_u \), arising respectively from alternating lifts of \( x^{m}(X)_P \) and \( y(X)_P \), and of \( x^{m}(X)_Q \) and \( z(X)_Q \). In Subsection~\ref{IOTZC}, the construction proceeded from \( C \) to \( \mathcal{U}_u \) and then to \( D_u \), whereas here we construct \( \mathcal{U}_u \) from \( C \) and \( D_u \). Regardless of the difference between the two construction orders, the resulting relative position of \( C \) and \( D_u \) is the same. As illustrated in Figure~\ref{fig_zigzag_def2}, we next define the points \( M_i, P_i', P_i'', H_i, N_i, Q_i', Q_i'' \) (for all integers \( i \)) and the geodesics \( \mathcal{U}_u \) and \( \mathcal{V}_u \).

Let \( M_0 \) be the intersection point of \( L \) with a segment in \( C \) that is a lift of some \( x^m \)-segment. Once \( M_0 \) is fixed, we can define the points \( M_i, P_i', P_i'', H_i \) for all integers \( i \) in the same manner as before (see Figure~\ref{fig_zigzag_def2}). Since \(\varepsilon_P(x,y) = -\varepsilon_Q(x,z)\), there exists a geodesic \(\mathcal{U}_u\) perpendicular to \( L \) such that \(\rho_{\mathcal{U}_u}(C) = D_u\), sending each \( x^m \)- and \( y \)-segment in \( C \) to the corresponding \( x^m \)- and \( z \)-segment in \( D_u \). Then we can define \(u, a, \theta_0, T', R', \bar{R}, \mathcal{V}_u'\) and \( N_i, Q_i', Q_i'' \) for all \( i \). In Figure~\ref{fig_zigzag_def2}, the segments corresponding to \( x^m \) (in place of \(\alpha\)) are drawn in blue, those corresponding to \( y \) (in place of \(\beta\)) are drawn in red, and those corresponding to \( z \) are drawn in green.

(i) Suppose \( \ell_{x^m} < \ell_y =\ell_z \) and \( \theta_{0} < \frac{\pi}{2} \).

\noindent
We classify the relative positions of the zigzag curves \( C \) and \( D_u \) into four cases: \(\mathrm{(I)}\) \( 0 < u < 2a \), \(\mathrm{(II)}\) \( u = 2a \), \(\mathrm{(III)}\) \( 2a < u < \frac{1}{2} \ell_{\left| x_{P}^m y_{P} \right|} \) , \(\mathrm{(IV)}\) \( u = \frac{1}{2} \ell_{\left| x_{P}^m y_{P} \right|} \) (see Figure \ref{fig_zigzag5_3}). By Lemma \ref{zigzag_lem2}, in Case \(\mathrm{(II)}\), \(\mathrm{(III)}\), and \(\mathrm{(IV)}\), there exist a \((P_0''P_1',\, Q_0'Q_0'')\)-transverse intersection point \( R'\) and a \((P_0'P_0'',\, Q_{-1}''Q_0')\)-transverse intersection point \(\bar{R}'\) such that \( \phi_{R'} = \phi_{\bar{R}'} < \phi_P \). Let \( R \coloneqq \mathrm{pr}(R') \) and \( \bar{R} \coloneqq \mathrm{pr}(\bar{R}')\). Then the \((x(X),\, y(X))\)-transverse intersection point \( R \) and the \((x(X),\, z(X))\)-transverse intersection point \( \bar{R} \) satisfy \( \phi_R = \phi_{\bar{R}} < \phi_P = \phi_Q \). Therefore, the claim holds in these cases.

In Case \(\mathrm{(I)}\), the conditions \( \theta_0 < \frac{\pi}{2} \) and \( 0 < u < 2a \) ensure that the geodesic \( \mathcal{U}_u \) intersects the segment \( P_0'P_0'' \) transversely. Moreover, since \( \rho_{\mathcal{U}_u}(C)^{-1} = D_u \), the segments \( P_0'P_0'' \) and \( Q_0'Q_0'' \) intersect transversely. We denote this intersection point by \( S' \). Here, the loop \( x(X) \) has a transverse self-intersection point \( S \coloneqq \mathrm{pr}(S') \). Since \( m \ge 2 \), several lifts of \( S \) lie on \( P_0'P_0'' \), and they appear at points whose distances from \( S' \) are positive integer multiples of \( \ell_x \). For each of these points, there exists a lift of \( x \) passing through it whose forward angle is \( \phi_S \). Let \( S'' \) be one of these lifts, distinct from \( S' \). Let \( L_x \) be the bi-infinite geodesic extension of the above lift of \( x \) passing through \( S'' \). If \( m \ge 3 \), we may choose \( S'' \) to lie in the interior of the segment \( P_0'P_0'' \). If \( m = 2 \), the position of \( S'' \) depends on the parameter \( u \). When \( u \neq a \), we may also choose \( S'' \) to lie in the interior of the segment \( P_0'P_0'' \). When \( u = a \), the lift \( S'' \) can be chosen only from the two endpoints \( P_0' \) or \( P_0'' \). Hence, we classify the position of \( S'' \) on the segment \( P_0'P_0'' \) into four cases: (I-A) \( S'' \) lies on the interior of the segment \( S'P_0'' \), (I-B) \( S'' \) lies on the interior of the segment \( S'P_0' \), (I-C) \( S'' \) lies on \( P_0'' \), (I-D) \( S'' \) lies on \( P_0' \).

In Case (I-A), the geodesic \( L_x \) enters the quadrilateral \( S'Q_0''T'P_0'' \) through the side \( S'P_0'' \). Since two geodesics in the upper half-plane intersect transversely at most once, \( L_x \) cannot exit the quadrilateral through the same side \( S'P_0'' \). If \( L_x \) exits through \( S'Q_0'' \), let \( S''' \) denote the intersection point. Since the sum of the interior angles of the triangle \( \triangle S'S''S''' \) exceeds \( \pi \), this contradicts the Gauss–Bonnet theorem. Hence, \( L_x \) must exit the quadrilateral through either \( T'Q_0'' \) or \( T'P_0'' \). We classify the exits of \( L_x \) from the quadrilateral into two cases: (I-A-1) \( L_x \) exits the quadrilateral from the segment \( T'Q_0'' \), (I-A-2) \( L_x \) exits the quadrilateral from the segment \( T'P_0'' \) (see Figure \ref{fig_zigzag_A}).

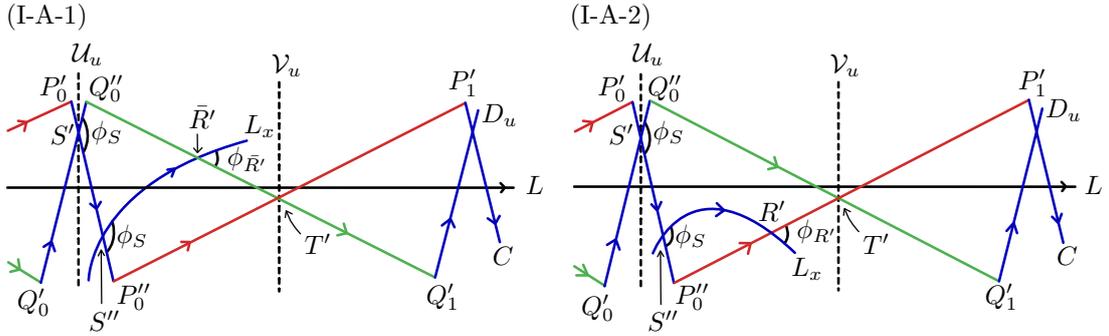
\begin{figure}[H]
\centering
\begin{overpic}[scale=0.97]{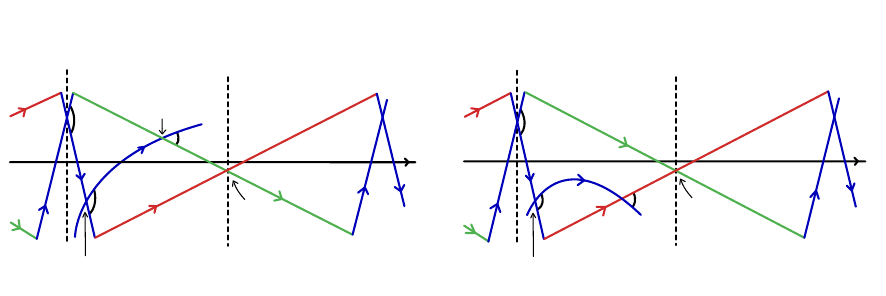}
\put(7,26.5){$\mathcal{U}_u$}
\put(25,25.5){$\mathcal{V}_u$}
\put(45,8){$C$}
\put(44,21){$D_u$}
\put(4,23.5){$P_0'$}
\put(41,24){$P_1'$}
\put(11,4.5){$P_0''$}
\put(28,9){$T'$}
\put(2,4){$Q_0'$}
\put(39,4.5){$Q_1'$}
\put(8.5,23.5){$Q_0''$}
\put(48,14.5){$L$}
\put(22.5,20){$L_x$}
\put(17.5,20.5){$\bar{R'}$}
\put(21,17.1){$\phi_{\bar{R'}}$}
\put(5,19){$S'$}
\put(8.5,2){$S''$}
\put(11,10){$\phi_{S}$}
\put(8.7,19.5){$\phi_{S}$}

\put(57.7,26.5){$\mathcal{U}_u$}
\put(75.5,25.5){$\mathcal{V}_u$}
\put(96,8){$C$}
\put(95,21.5){$D_u$}
\put(54.5,23.5){$P_0'$}
\put(92.5,24){$P_1'$}
\put(61.5,4.5){$P_0''$}
\put(78.5,9){$T'$}
\put(53,3.5){$Q_0'$}
\put(89.5,4.5){$Q_1'$}
\put(59,23.5){$Q_0''$}
\put(98.5,14.5){$L$}
\put(72,7){$L_x$}
\put(55.5,19){$S'$}
\put(59,2){$S''$}
\put(61.5,10){$\phi_{S}$}
\put(59.5,19.5){$\phi_{S}$}
\put(72.5,10.5){$\phi_{R'}$}
\put(69,12){$R'$}

\put(1,30){(I-A-1)}
\put(52,30){(I-A-2)}
\end{overpic}
\caption{The case where the geodesic \(L_x\) exits from the segment \(T'Q_{0}''\) (left) and the case where it exits from the segment \(T'P_{0}''\) (right).}
\label{fig_zigzag_A}
\end{figure}

In Case (I-A-1), \( L_x \) exits the quadrilateral from \( T'Q_0'' \), and we denote the intersection point by \( \bar{R}' \). Then \( \bar{R} \coloneqq \mathrm{pr}(\bar{R}') \) is the \((x(X), z(X))\)-transverse intersection point. In the quadrilateral \( Q_0''S'S''\bar{R}' \), the sum of the interior angles is less than \( 2\pi \), and we obtain \( \phi_{\bar{R}} < \phi_P = \phi_Q \). Moreover, since \( \rho_{\mathcal{U}_u}(C) = D_u^{-1} \), there exists a corresponding \((x(X), y(X))\)-transverse intersection point \( R \) satisfying \( \phi_R = \phi_{\bar{R}} \).

In Case (I-A-2), \( L_x \) exits the quadrilateral from \( T'P_0'' \), and we denote this intersection point by \( R' \). Then \( R \coloneqq \mathrm{pr}(R') \) is the \((x(X), y(X))\)-transverse intersection point. In the triangle \( \triangle S''P_0''R' \), the sum of the interior angles is less than \( \pi \), and we obtain \(\phi_R < \phi_P = \phi_Q\).
Again, since \( \rho_{\mathcal{U}_u}(C) = D_u^{-1} \), there exists a corresponding \((x(X), z(X))\)-transverse intersection point \( \bar{R} \) satisfying \( \phi_{\bar{R}} = \phi_R \).

In Case (I-B), we can apply the same argument to the quadrilateral \( \rho_{\mathcal{U}_u}(S'Q_0''T'P_0'') \) as in Case (I-A).

In Case (I-C), the point \( P_0'' \) lies on three geodesics: \( P_0'P_0'' \), \( P_0''P_1' \), and \( L_x \). We classify this case into three cases by comparing \( \phi_S \) and \( \phi_P \): (I-C-1) \( \phi_S < \phi_P \), (I-C-2) \( \phi_S = \phi_P \), (I-C-3) \( \phi_S > \phi_P \).

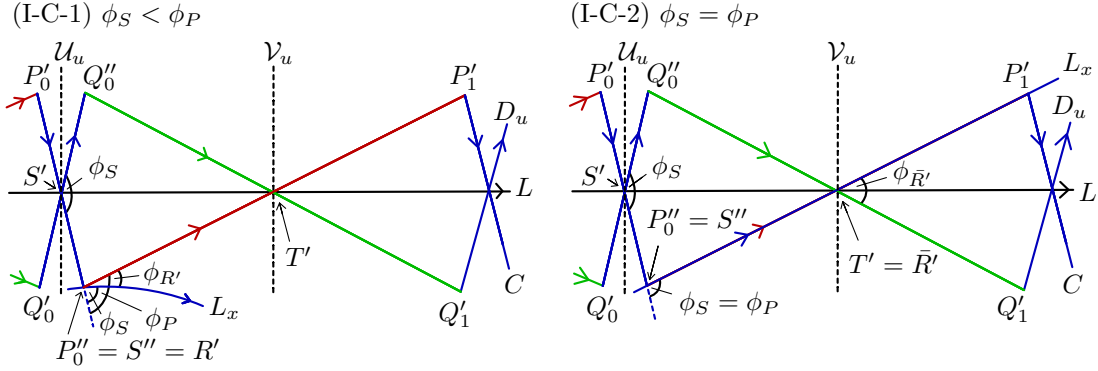
\begin{figure}[H]
\centering
\begin{overpic}[scale=0.4]{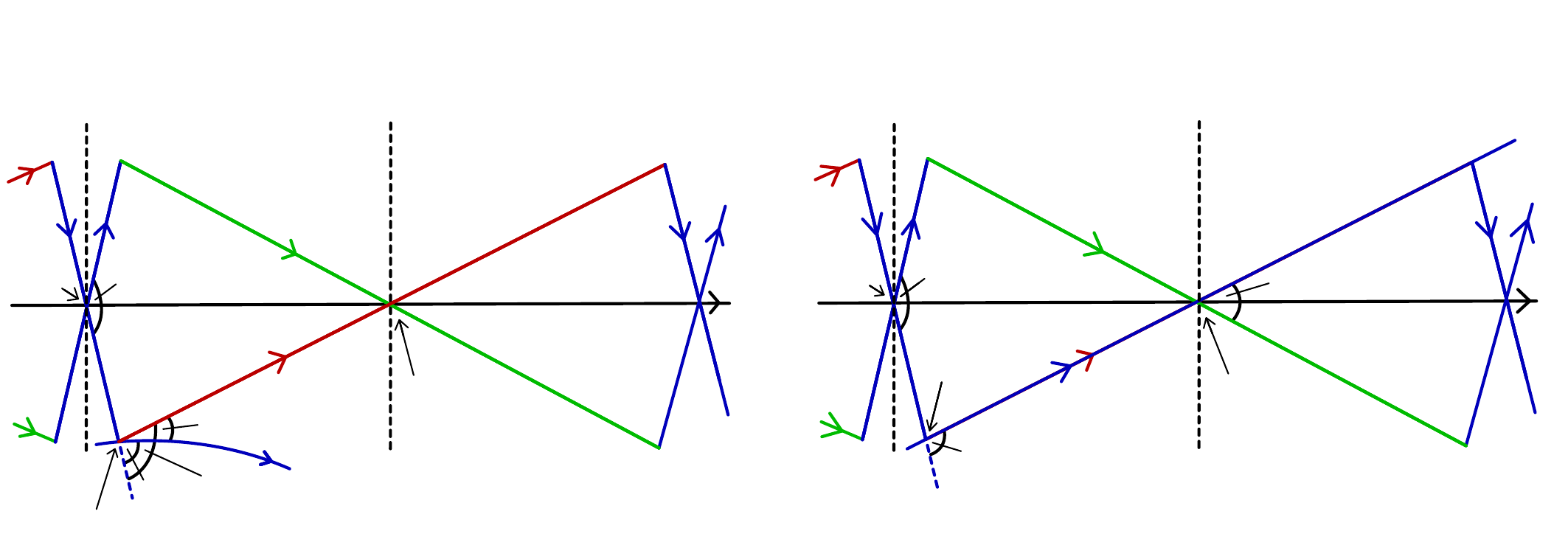}
\put(5,27.5){$\mathcal{U}_u$}
\put(24,27.5){$\mathcal{V}_u$}
\put(46,6){$C$}
\put(45,22){$D_u$}
\put(2,25.3){$P_0'$}
\put(41,25.3){$P_1'$}
\put(5,0){$P_0''=S''=R'$}
\put(26,8.5){$T'$}
\put(2,4){$Q_0'$}
\put(40,3.5){$Q_1'$}
\put(7.5,25.3){$Q_0''$}
\put(47,15){$L$}
\put(19,4){$L_x$}
\put(2,16){$S'$}
\put(9,2.6){$\phi_{S}$}
\put(13.2,3){$\phi_{P}$}
\put(8,17){$\phi_{S}$}
\put(13,7.2){$\phi_{R'}$}

\put(56.5,27.5){$\mathcal{U}_u$}
\put(75.5,27.5){$\mathcal{V}_u$}
\put(97,6){$C$}
\put(96,22.5){$D_u$}
\put(53.5,25.5){$P_0'$}
\put(91.5,25.5){$P_1'$}
\put(59,12){$P_0''=S''$}
\put(77.5,7.9){$T'=\bar{R'}$}
\put(53.5,4){$Q_0'$}
\put(91,4){$Q_1'$}
\put(59,25.5){$Q_0''$}
\put(98.5,14.5){$L$}
\put(97,26.5){$L_x$}
\put(53,16){$S'$}
\put(62,4.5){$\phi_{S}=\phi_{P}$}
\put(59.5,17){$\phi_{S}$}
\put(81.5,16.5){$\phi_{\bar{R'}}$}

\put(1,31){(I-C-1) \( \phi_S < \phi_P \)}
\put(52,31){(I-C-2) \( \phi_S = \phi_P \)}
\end{overpic}
\caption{The case where \(S''\) lies on \(P_0''\) with \(\phi_S < \phi_P\) (left) and the case where \(S''\) lies on \(P_0''\) with \(\phi_S = \phi_P\) (right).}
\label{fig_zigzag_C}
\end{figure}

In Case (I-C-1), the geodesic \( L_x \) intersects \( P_0''P_1' \) transversely, and we denote this intersection point by \( R' \) (See Figure \ref{fig_zigzag_C}). Then \( R \coloneqq \mathrm{pr}(R') \) is the \((x(X), y(X))\)-transverse intersection point. We have \( \phi_R = \phi_P - \phi_S < \phi_P \). Moreover, since \( \rho_{\mathcal{U}_u}(C) = D_u^{-1} \), there exists a corresponding \((x(X), z(X))\)-transverse intersection point \( \bar{R} \) satisfying \( \phi_{\bar{R}} = \phi_R \).

In Case (I-C-2), the geodesic \( L_x \) contains the segment \( P_0'P_0'' \), and hence intersects \( Q_0'Q_0'' \) transversely (See Figure \ref{fig_zigzag_C}). We denote this intersection point by \( \bar{R}' \). Then \( \bar{R} \coloneqq \mathrm{pr}(\bar{R}') \) is the \((x(X), y(X))\)-transverse intersection point. In the quadrilateral \( Q_0''S'S''\bar{R}' \), the sum of the interior angles is less than \( 2\pi \). Combining this inequality with the assumption \( \phi_S = \phi_P \), we obtain \( \phi_{\bar{R}} < \phi_P = \phi_Q \). Moreover, since \( \rho_{\mathcal{U}_u}(C) = D_u^{-1} \), there exists a corresponding \((x(X), y(X))\)-transverse intersection point \( R \) satisfying \( \phi_R = \phi_{\bar{R}} \).

In Case (I-C-3), the geodesic \( L_x \) enters into the quadrilateral \( Q_0''S'S''T' \) and intersects \( P_0''P_1' \) transversely at the point \( P_0'' \). Hence, \( L_x \) exits the quadrilateral from the segment \( T'Q_0'' \). Then the same argument as in Case (I-C-1) applies.

In Case (I-D), the points \( P_0'' \) and \( P_0' \) both lie on the lift \( P_0'P_0'' \) of \( x^m(X) \), and the distance between them is \( m\ell_x \), which is a positive integer multiple of \( \ell_x \). Hence there exists another lift of \( x \) that passes through \( P_0'' \) with the forward angle \( \phi_S \). Then the same argument as in Case (I-C) applies.

(ii) Suppose \( \ell_{x^m} < \ell_y =\ell_z \) and \( \theta_0 = \frac{\pi}{2} \).

\noindent
By Lemma \ref{zigzag_lem3}, in Case \(\mathrm{(V)}\) and \(\mathrm{(VI)}\), there exist a \((P_0''P_1',\, Q_0'Q_0'')\)-transverse intersection point \( R'\) and a \((P_0'P_0'',\, Q_{-1}''Q_0')\)-transverse intersection point \(\bar{R}'\) such that \( \phi_{R'} = \phi_{\bar{R}'} < \phi_P \). Let \( R \coloneqq \mathrm{pr}(R') \) and \( \bar{R} \coloneqq \mathrm{pr}(\bar{R}') \). Then the \((x(X),\, y(X))\)-transverse intersection point \( R \) and the \((x(X),\, z(X))\)-transverse intersection point \( \bar{R} \) satisfy \( \phi_R = \phi_{\bar{R}} < \phi_P = \phi_Q \). Hence, the claim holds in these cases.

(iii) Suppose \(  \ell_{x^m} < \ell_y =\ell_z \) and \( \theta_0 > \frac{\pi}{2} \).

\noindent
In this case, we classify the relative position of the zigzag curves \( C \) and \( D_u \) into a single case: \(\mathrm{(VII)}\) (See Figure \ref{fig_zigzag7}). By Lemma \ref{zigzag_lem3}, in Case \(\mathrm{(VII)}\), there exist a \((P_0''P_1',\, Q_0'Q_0'')\)-transverse intersection point \(R'\) and a \((P_0'P_0'',\, Q_{-1}''Q_0')\)-transverse intersection point \(\bar{R}'\) such that \(\phi_{R'} = \phi_{\bar{R}'} < \phi_P\). Let \(R \coloneqq \mathrm{pr}(R')\) and \(\bar{R} \coloneqq \mathrm{pr}(\bar{R}')\). Then the \((x(X), y(X))\) and \((x(X), z(X))\)-transverse intersection points \(R\) and \(\bar{R}\) satisfy \(\phi_R = \phi_{\bar{R}} < \phi_P = \phi_Q\). Therefore, the claim holds in this case.

(iv) Suppose \( \ell_{x^m} = \ell_y =\ell_z \).

\noindent
In this case, we classify the relative positions of the zigzag curves \( C \) and \( D_u \) into two cases: \(\mathrm{(VIII)}\) \( 0 < u < \frac{1}{2} \ell_{\left| x_{P}^m y_{P} \right|} \), \(\mathrm{(IX)}\)  \( u = \frac{1}{2} \ell_{\left| x_{P}^m y_{P} \right|} \) (See Figure \ref{fig_zigzag8}).

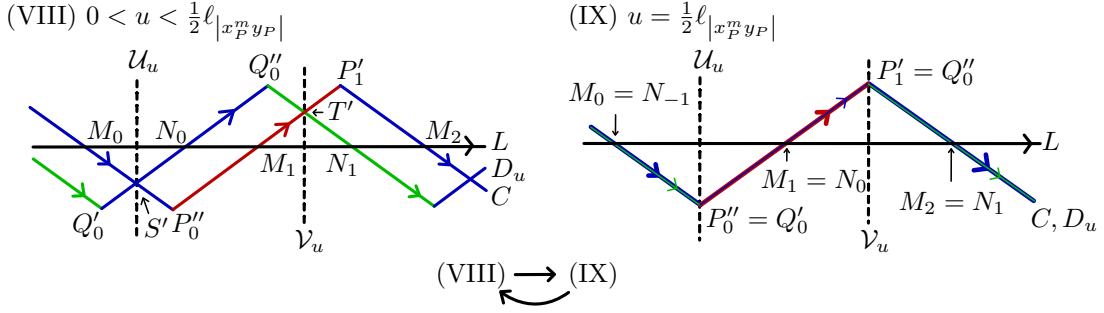
\begin{figure}[H]
\centering
\begin{overpic}[scale=0.58]{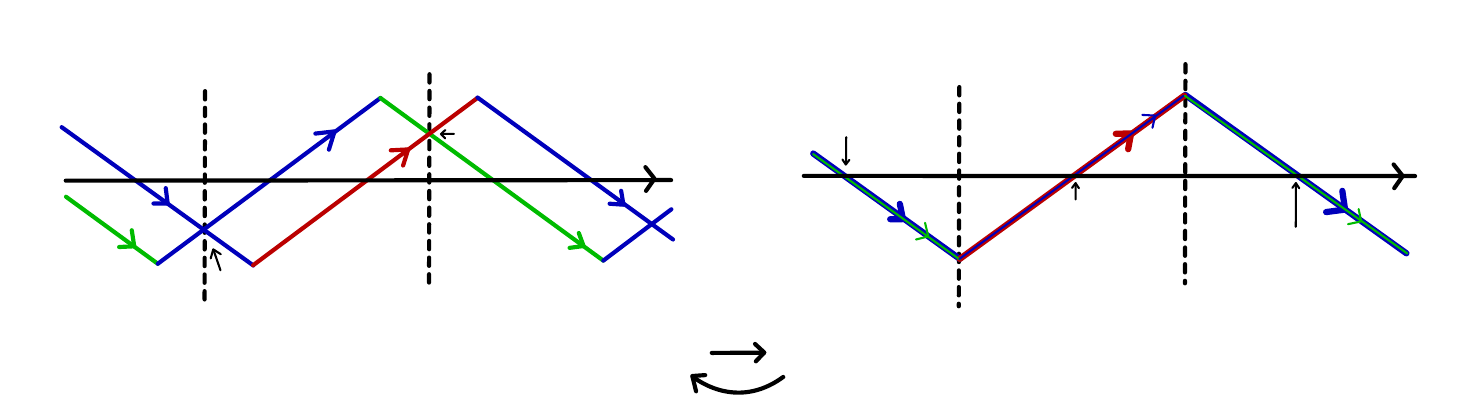}
\put(2,27){\(\mathrm{(VIII)}\) \( 0 < u < \frac{1}{2} \ell_{\left| x_{P}^m y_{P} \right|} \)}
\put(53,27){\(\mathrm{(IX)}\)  \( u = \frac{1}{2} \ell_{\left| x_{P}^m y_{P} \right|} \)}

\put(13.2,23){$\mathcal{U}_u$}
\put(28.2,6.5){$\mathcal{V}_u$}
\put(31.2,18.2){$T'$}
\put(14.5,7.5){$S'$}
\put(17,8){$P_{0}''$}
\put(8,8){$Q_{0}'$}
\put(32,22){$P_{1}'$}
\put(24,22.5){$Q_{0}''$}
\put(46,10.5){$C$}
\put(46,13){$D_u$}
\put(46,15.5){$L$}
\put(15.5,16.5){$N_0$} \put(9.3,16.5){$M_0$} \put(25,13.5){$M_1$} \put(31,13.5){$N_1$} \put(40,16.5){$M_2$}

\put(64.5,23){$\mathcal{U}_u$}
\put(79.5,6.5){$\mathcal{V}_u$}
\put(95,8.5){$C,D_u$}
\put(96,15.5){$L$}
\put(65.5,8.5){$P_0''=Q_0'$}
\put(81,22){$P_1'=Q_0''$}
\put(53,20){$M_0=N_{-1}$}
\put(70.5,12){$M_1=N_{0}$}
\put(83,10){$M_2=N_{1}$}

\put(41,3.4){$\mathrm{(VIII)}$}
\put(53,3.4){$\mathrm{(IX)}$}
\end{overpic}
\caption{Possible configurations of the zigzag curves $C$ and $D_u$ in the case \( \varepsilon_P(x, y) = -\varepsilon_Q(x, z) \), \( \ell_y = \ell_z = \ell_{x^m} \)}
\label{fig_zigzag8}
\end{figure}

In Case \(\mathrm{(VIII)}\), by the assumption \( \ell_{x^m} = \ell_y = \ell_z \), for each integer \( i \), the triangle \( \triangle M_i P_i' M_{i+1} \) satisfies \( M_i P_i' = M_{i+1} P_i' = \frac{1}{2}\,\ell_{\left| x_{P}^m y_{P} \right|} \). Since \( H_i \) is the midpoint of \( M_i M_{i+1} \), we obtain \( 2a = \frac{1}{2}\,\ell_{\left| x_{P}^m y_{P} \right|} \).
Hence, the assumption \( 0 < u < \frac{1}{2}\,\ell_{\left| x_{P}^m y_{P} \right|} \)  
is the condition \( 0 < u < 2a \) which ensures that the geodesic \( \mathcal{U}_u \) intersects the segment \( P_0'P_0'' \) transversely. Moreover, since \( \rho_{\mathcal{U}_u}(C)^{-1} = D_u \), the segments \( P_0'P_0'' \) and \( Q_0'Q_0'' \) intersect transversely. We denote this intersection point by \( S' \). Then the same argument as in Case \(\mathrm{(I)}\) applies.

In Case \(\mathrm{(IX)}\), we obtain \( y = z = x^m \).

Finally, we consider the case \( \ell_y = \ell_z < \ell_{x^m} \). However, due to the argument given in the case \( \ell_y = \ell_z > \ell_{x^m} \), we hardly need to discuss this case separately. Indeed, it suffices to switch the red and blue segments of the zigzag curve \( C \), and similarly to switch the green and blue segments of the zigzag curve \( D_u \). We then denote the corresponding cases to \(\mathrm{(I)}\) through \(\mathrm{(VII)}\) by \(\mathrm{(I)}'\) through \(\mathrm{(VII)}'\), respectively. In particular, by observing that the forward angles at the points \( R' \) and \( \bar{R}' \) remain unchanged under this transformation, the claims in all cases \(\mathrm{(II)}'\) through \(\mathrm{(VII)}'\) follow immediately. In Case \(\mathrm{(I)}'\), we simply exchange the roles of \(T'\) and \(S'\) in Case \(\mathrm{(I)}\). Then the same argument as in Case \(\mathrm{(I)}\) applies and the proof is complete.
\end{proof}

\begin{lemma} \label{xyz_+}
Let \( X \) be a complete hyperbolic metric, and let \( x, y, z \in \hat{\pi} \) be essential loops with \( \ell_y(X) = \ell_z(X) \). If there exist an \((x(X),\, y(X))\)-transverse intersection point \(P\) and an \((x(X),\, z(X))\)-transverse intersection point \(Q\) such that
\( \varepsilon_P(x, y) = \varepsilon_Q(x, z) \) and \( \left| x_P y_P \right| = \left| x_Q z_Q \right| \), then there exist an \((x(X),\, y(X))\)-transverse intersection point \(R\) and an \((x(X),\, z(X))\)-transverse intersection point \(\bar{R}\) such that
\[
\phi_R = \phi_{\bar{R}} < \phi_P = \phi_Q, 
\quad \text{or} \quad y = z.
\]
\end{lemma}

\begin{proof}
As in the proof of Lemma \ref{xyz_-}, we obtain zigzag curves \( C \) and \( D \) along a common bi-infinite geodesic line \( L \subset \mathbb{H} \), now with \(x^m\) replaced by \(x\).
Fix a point \( M_0 \) on an \(x\)-segment of \( C \) that meets \( L \), and define \( M_i, P_i', P_i'' \) as before (see Figure~\ref{fig_zigzag0}).

\begin{figure}[H]
\centering
\begin{overpic}[scale=0.6]{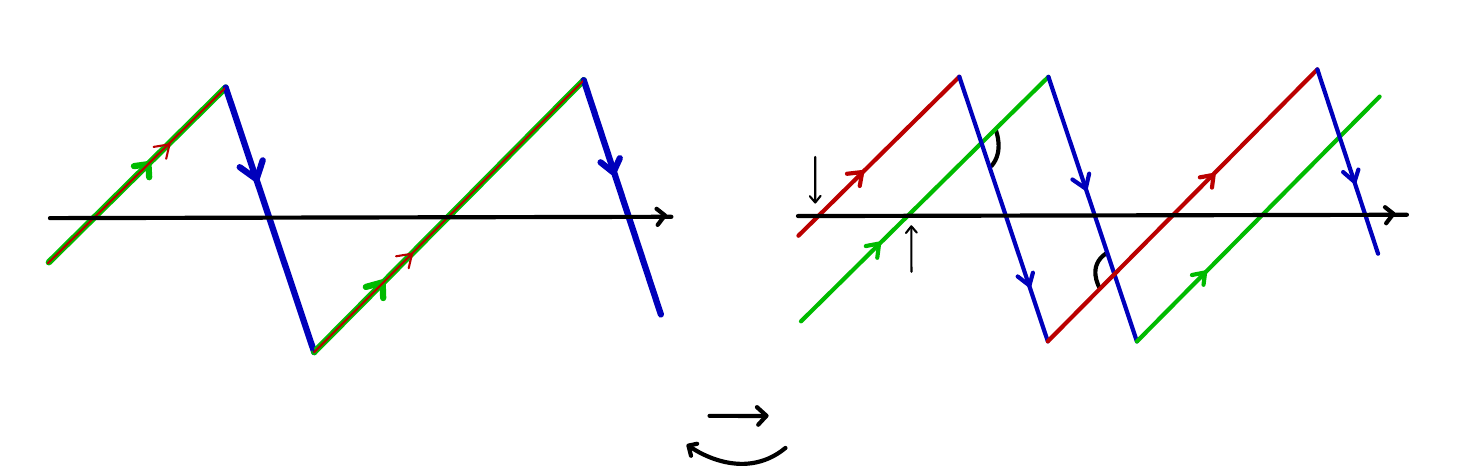}
\put(2,30){$\mathrm{(X)}$ $P_0'=Q_0'$}
\put(53,30){$\mathrm{(XI)}$ $P_0'\neq Q_0'$}

\put(16,26){$P_0'=Q_0'$}
\put(21,6){$P_0''=Q_0''$}
\put(40.5,26){$P_1'=Q_1'$}
\put(46,17){$L$}
\put(42.5,8){$C,D$}
\put(5.6,14.5){$M_{-1}=N_{-1}$} \put(18.5,18){$M_0=N_{0}$} \put(29.6,14.5){$M_{1}=N_{1}$} \put(34,11){$M_2=N_{2}$}

\put(65.5,26.5){$P_0'$}
\put(71.5,26.5){$Q_0'$}
\put(70,6){$P_0''$}
\put(75,6){$Q_0''$}
\put(89.5,27){$P_1'$}
\put(76.3,12.3){$R'$}
\put(63,22){$\bar{R}'$}
\put(70.6,13.7){$\phi_{R'}$}
\put(68,21.5){$\phi_{\bar{R}'}$}
\put(93,13){$C$}
\put(93,26){$D$}
\put(96,17){$L$}
\put(54,22){$M_{-1}$} \put(60,11){$N_{-1}$} \put(74.5,18.2){$N_0$} \put(65,15){$M_0$} \put(79,15){$M_1$} \put(85,15){$N_1$} \put(89.5,15){$M_2$}

\put(43,3.4){$\mathrm{(X)}$}
\put(53,3.4){$\mathrm{(XI)}$}
\end{overpic}
\caption{Possible configurations of the zigzag curves $C$ and $D_u$ in the case \( \varepsilon_P(x, y) = \varepsilon_Q(x, z) \)}
\label{fig_zigzag2}
\end{figure}
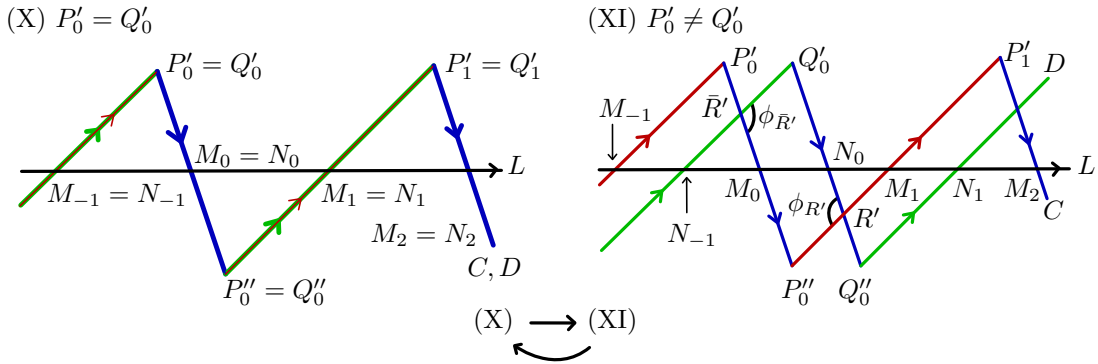

By periodicity, we may fix \( C \) and shift \( D \) along \( L \) until it is displaced by the distance \( \ell_{\left| x_{P} y_{P} \right|} \).

We first place \( D \) so that \( Q_0' \) coincides with \( P_0' \). At this point, \( y = z \) (see Case $\mathrm{(X)}$ in Figure \ref{fig_zigzag2}). Then, moving \( D \) along \( L \) by a distance \( d \) with \( 0 < d < \ell_{\left| x_{P}^m y_{P} \right|} \), the segment \( Q_0'Q_0'' \) intersects \( P_1'P_0'' \)(see Case $\mathrm{(XI)}$ in Figure \ref{fig_zigzag2}). Since \( Q_0'Q_0'' \) enters either \( \triangle M_0P_0''M_1 \) or \( \triangle M_1P_1'M_2 \), it must exit from one of these triangles. By the Gauss–Bonnet theorem, it must exit through \( P_1'P_0'' \). Let this intersection point be \( {R}' \). Similarly, the segment \( P_0'P_0'' \) intersects \( Q_0'Q_{-1}'' \), and denote the intersection point by \( \bar{R}' \). Let \( R \coloneqq \mathrm{pr}(R') \). Then \( R \) is the \((x(X),y(X))\)-transverse intersection point, and \( R' \coloneqq \mathrm{pr}(\bar{R}') \) is the \((x(X),z(X))\)-transverse intersection point. By symmetry, we have \( \phi_R = \phi_{\bar{R}} \).
\end{proof}

By substituting $y$ for $z$ in Lemma \ref{xyz_-}, we obtain the following lemma.

\begin{lemma}\label{xy}
Let \( X \) be a complete hyperbolic metric. Let \(x, y \in \hat{\pi}\) be essential loops, and suppose that \(P\) and \(Q\) are \((x(X), y(X))\)-transverse intersection points satisfying \(\varepsilon_{P}(x,y) = -\varepsilon_{Q}(x,y)\). If there exists a positive integer \( m \geq 2 \) such that \( \left| x^{m}_{P} y_{P} \right| = \left| x^{m}_{Q} y_{Q} \right| \),
then either there exists a point \( R \in x(X) \cap y(X) \) such that
\[
\phi_R < \phi_P = \phi_Q \quad \text{or} \quad y = x^m.
\]
\end{lemma}

The next lemma follows from Lemma \ref{xyz_-} and \ref{xyz_+}

\begin{lemma} \label{xyz}
Let \( X \) be a complete hyperbolic metric, and \( x, y, z \in \hat{\pi} \) with \( \ell_y(X) = \ell_z(X) \). Suppose that \(P\) is an \((x(X), y(X))\)-transverse intersection point and \(Q\) is an \((x(X), z(X))\)-transverse intersection point. If there exists a positive integer \( m \geq 2\) such that \( \left| x_P^m y_P \right| = \left| x_Q^m z_Q \right| \), then there exist an \((x(X), y(X))\)-transverse intersection point \(R\) and an \((x(X), z(X))\)-transverse intersection point \(\bar{R}\) such that
\[
\phi_R = \phi_{\bar{R}} < \phi_P = \phi_Q, 
\quad \text{or} \quad y = z.
\]
\end{lemma}

\subsection{Proof of the separability criteria}
Now we prove the separability criteria, Theorems~\ref{WSC1} and~\ref{SSC1}. To prove these theorems, we begin with the following two lemmas. The proof of the first one is based on Lemma~\ref{xy}.

\begin{lemma} \label{KeyLemma1}
Let $x, y \in \hat{\pi}$ be essential loops. If there exists $m \in \mathbb{N}_{\geq 2}$ such that $[x^m, y] = 0 \in K\hat{\pi}$, then either $y = x^m$ or there is no $(x(X), y(X))$-transverse intersection point for any complete hyperbolic metric $X$.
\end{lemma}

\begin{proof}
Assume \( y\neq x^m \) and there is an $(x(X), y(X))$-transverse intersection point for a complete hyperbolic metric $X$. By the definition of the Goldman bracket,
\[
[x^m, y] = \sum_{P \in x^m \cap y} \varepsilon_{P}(x^m,y)|{x^m_{P}y_{P}}| = m \sum_{P \in x \cap y} \varepsilon_{P}(x,y)|{x^m_{P}y_{P}}| .
\]

Here, the representatives of \( x^m \) and \( y\) are taken to be the geodesic representatives with respect to the complete hyperbolic metric \( X \). Choose an $(x(X), y(X))$-transverse intersection point $P$ whose forward angle $\phi_P$ is minimum among all $(x(X), y(X))$-transverse intersection points, and the corresponding term $|{x^m_{P}y_{P}}|$. Since we have \( [x^m, y] = 0 \),
\[
|{x^m_{P}y_{P}}|=|{x^m_{Q}y_{Q}}|, \quad \varepsilon_{P}(x,y)=-\varepsilon_{Q}(x,y).
\]
for some $(x(X), y(X))$-transverse intersection point $Q$. By Lemma \ref{xy}, there exists an $(x(X), y(X))$-transverse intersection point $R$ such that the forward angle $x$ and $y$ at $R$ is strictly smaller than that at $P$, i.e., $\phi_R < \phi_P$. This contradicts the assumption that $P$ was minimum among all \((x(X), y(X))\)-transverse intersection points. This completes the proof.
\end{proof}

\begin{lemma} \label{KeyLemma2}
Let $X$ be a complete hyperbolic metric, and let $x, y \in \hat{\pi}$ be essential loops. If there is no $(x(X), y(X))$-transverse intersection point, then one of the following conditions \textnormal{(A)} or \textnormal{(B)} holds.\\
{\rm (A)} $x(X) \cap y(X) = \emptyset$,\\
{\rm (B)} There exists a simple and primitive loop $z \in \hat{\pi}$ and integers $m_1, m_2 \in \mathbb{Z}\setminus \{0\}$ such that
\[
x(X) = z^{m_1}(X), \quad y(X) = z^{m_2}(X).
\]
In both cases, the geometric intersection number satisfies $i(x, y) = 0$.
\end{lemma}

\begin{proof}
Assume that \textnormal{(A)} does not hold.
Then $x(X) \cap y(X) \neq \emptyset$.
Since $x(X)$ and $y(X)$ are both closed geodesics, there exists a primitive element $z \in \hat{\pi}$ and integers $m_1, m_2 \in \mathbb{Z} \setminus \{0\}$ such that
\[
x(X) = z^{m_1}(X), \quad y(X) = z^{m_2}(X).
\]

Moreover, assume that $z(X)$ has a transverse self-intersection point. Then $x(X)$ and $y(X)$ have an $(x(X), y(X))$-transverse intersection point, contradicting our assumption. Hence, $z(X)$ is a simple closed curve.

By the definition of the geometric intersection number, we have $i(x, y) = 0$ in case (A). In case (B), by taking a sufficiently small tubular neighborhood $N$ of $z(X)$, we can deform $x(X)$ and $y(X)$ within $N$ so that they lie on disjoint parallel copies of $z(X)$, showing again that $i(x, y) = 0$.
\end{proof}
By Lemma \ref{KeyLemma1} and Lemma \ref{KeyLemma2}, we obtain the following theorem, which is our first criterion.

\begin{thm}[Theorem \ref{WSC1}] \label{WSC2}
Let \( x, y \in \hat{\pi} \) and \( m \in \mathbb{N}_{\geq 2} \).  
Then \( [x^m, y] = 0 \) in \( K\hat{\pi} \) if and only if \( i(x, y) = 0 \) or \( y = x^m \).
\end{thm}

\begin{proof}
If $i(x, y) = 0$ or $y = x^m$, then $[x^m, y] = 0$ follows directly from the definition of the Goldman bracket. For the converse direction, the claim holds when either $x$ or $y$ is non-essential since $i(x, y) = 0$ in this case. If $x$ and $y$ are essential, then Lemma~\ref{KeyLemma1} and Lemma~\ref{KeyLemma2} imply that $i(x, y) = 0$ or $y = x^m$. This completes the proof.
\end{proof}

As a corollary of Theorem~\ref{WSC2}, we recover the following result of Chas and Kabiraj~\cite{Chas-Kabiraj2023}.

\begin{cor}\cite{Chas-Kabiraj2023} \label{1_m}
A primitive class $x \in \hat{\pi}$ contains a simple representative if and only if $[x, x^m] = 0$ for some \( m \ge 2 \).
\end{cor}

In the following example, for a free homotopy class $x \in \hat{\pi}$ and for every natural number $m \geq 2$, we verify that the Goldman bracket $[x^m, x]$ does not vanish.

\begin{ex}
Let $\Sigma$ be an oriented pair of pants as in Figure~\ref{SC16}, and \(a, b \in \pi_1(\Sigma, *)\) generators taken as shown in the same figure. As shown in Figure~\ref{SC16}, we take two representatives of $x \in \hat{\pi}$. Then, we choose one of them and, by going around it \(m\) times, obtain a representative of \(x^m\) for each natural number \(m \geq 2\). By using the method introduced in the Appendix of Chas and Krongold~\cite{Chas-Krongold2010}, we have
\[
[x^m, x] = m \bigl((b^{-1}a)^m ab^{-1} - (ab^{-1})^m b^{-1}a \bigr).
\]
Here, each term on the right-hand side represents an element of the set of all conjugacy classes 
\(\pi_1(\Sigma, *) / \mathrm{conj}\), and the two intersection points \(P\) and \(Q\) are each counted with multiplicity \(m\).

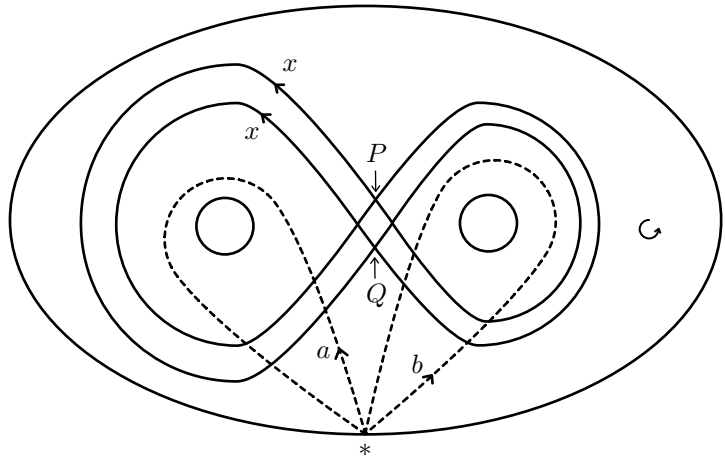
\begin{figure}[H]
\centering
\begin{overpic}[scale=0.5]{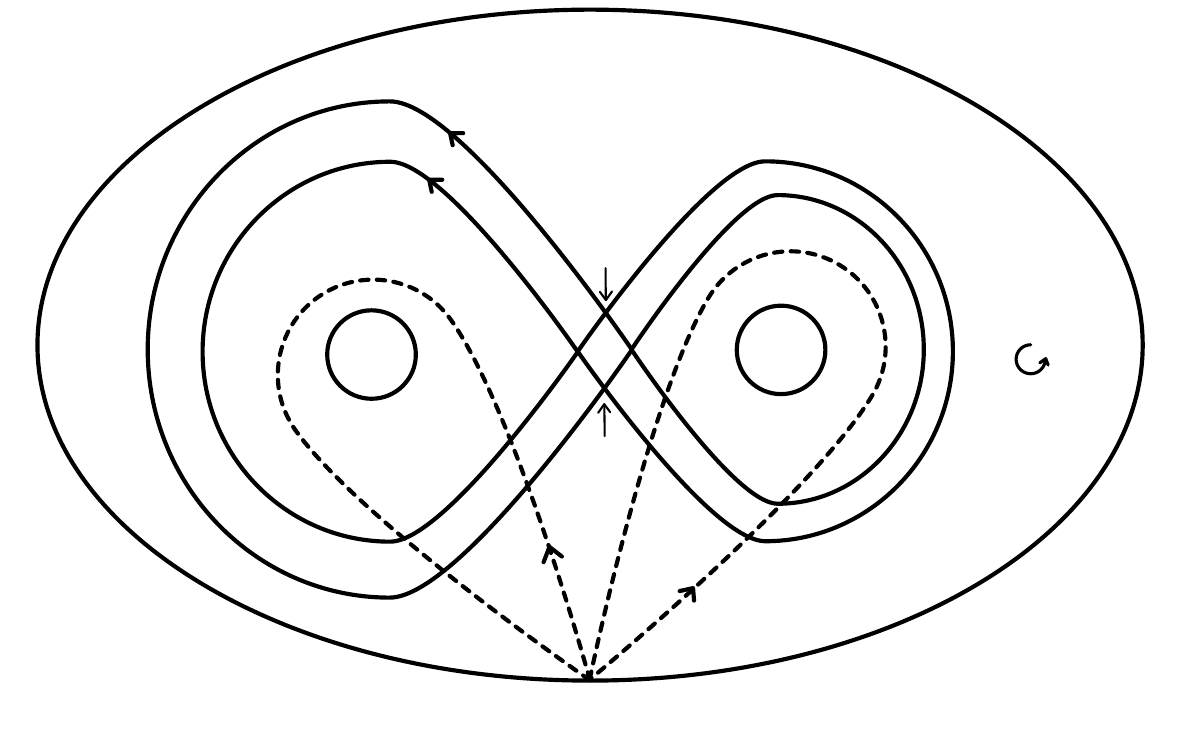}
\put(43.5,15){$a$}
\put(56,13){$b$}
\put(34,44){$x$}
\put(39,53){$x$}
\put(50,41){$P$}
\put(50,22.5){$Q$}
\put(49,2){$*$}
\end{overpic}
\caption{Two representatives of the free homotopy class $x \in \hat{\pi}$ on an oriented pair of pants}
\label{SC16}
\end{figure}

Let $w_1 = (b^{-1}a)^m ab^{-1} \in \pi_1(\Sigma,*)/\mathrm{conj}$ and 
$w_2 = (ab^{-1})^m b^{-1}a \in \pi_1(\Sigma,*)/\mathrm{conj}$.  
Both of them are cyclically reduced words, and each contains exactly one occurrence of $a^2$.  
By rewriting them so that they start with $a^2$, we have
\[
w_1 = a^2 b^{-2} (ab^{-1})^{m-1}, \qquad 
w_2 = a^2 (ab^{-1})^{m-1} b^{-2}.
\]
Since $m \geq 2$, the third letter from the left is $b^{-1}$ in $w_1$ but $a$ in $w_2$, so we have $w_1 \neq w_2$. Therefore $[x^m, x] \neq 0$.
\end{ex}

Since any Lie bracket is skew-symmetric, the vanishing of $[x^m, y]$ does not imply that $x$ and $y$ are disjoint. In fact, we have $[x^m, y]=0$ for $y=x^{m}$. Theorem~\ref{WSC2} states that apart from this trivial vanishing case of the bracket, the vanishing of the bracket does imply separability. In this sense, Theorem~\ref{WSC2} provides a ``weak'' separability criterion.

To fully determine whether two loops are separable, it suffices to verify any of the equivalent conditions in the next theorem. In this sense, Theorem~\ref{SSC2} provides ``strong" separability criteria. Before stating the theorem, we recall a fact that will be used in its proof. Since the fundamental group of \(\Sigma\) is torsion-free, the covering space corresponding to a nontrivial \(x \in \hat{\pi}\) is an annulus. Hence \(x^m = x^n \in \hat{\pi}\) implies \(m = n\) if \(x\) is nontrivial.

\begin{thm}[Theorem \ref{SSC1}] \label{SSC2}
Let $x, y \in \hat{\pi}$. The following four conditions are equivalent,\\
{\rm (1)} $i(x, y) = 0$,\\
{\rm (2)} There exist distinct $m_1, m_2 \in \mathbb{N}_{\geq 1}$ such that $[x^{m_1}, y] = [x^{m_2}, y] = 0 $ in \( K\hat{\pi} \),\\
{\rm (3)} There exist distinct $m_1, m_2 \in \mathbb{N}_{\geq 1}$ such that $[x^{m_1}, y] = [x, y^{m_2}] = 0$ in \( K\hat{\pi} \),\\
{\rm (4)} There exists $m \in \mathbb{N}_{\geq 2}$ and non-zero divisors $c_1, c_2 \in K$ such that $[x^m, c_1 y + c_2 y^{-1}] = 0$ in \( K\hat{\pi} \).
\end{thm}

\begin{proof}
If condition~(1) holds, then conditions~(2), (3), and (4) follow immediately. We now prove the converses. In proving the converses, it suffices to consider the case where both \(x\) and \(y\) are essential, as condition~(1) is satisfied when one of them is non-essential.

Assume that condition~(2) holds.  
We may assume \( m_1 > m_2 \).  
By Theorem~\ref{WSC2}, we have either \( y = x^{m_1} \) or \( i(x, y) = 0 \).  
If \( y = x^{m_1} \), we have
\[
[x^{m_2}, x^{m_1}] = 0.
\]
Then, by Lemma~\ref{KeyLemma1}, we have either \( x^{m_2} = x^{m_1} \)   
or there is no $(x^{m_2}(X), x^{m_1}(X))$-intersection point for any hyperbolic metric $X$.  
Since \( m_1 \ne m_2 \), we have \( x^{m_2} \ne x^{m_1} \).  
From the latter condition, it also follows that there is no $(x(X), y(X))$-transverse intersection point for any hyperbolic metric $X$.  
By Lemma~\ref{KeyLemma2}, we obtain $i(x, y) = 0.$

Assume that condition~(3) holds.  
We may again assume \( m_1 > m_2 \).  
By Theorem~\ref{WSC2}, we have either \( y = x^{m_1} \) or \( i(x, y) = 0 \).  
If \( y = x^{m_1} \), we have
\[
[x, x^{m_1 m_2}] = 0.
\]
Applying Lemma~\ref{KeyLemma1} again, we obtain either \( x = x^{m_1 m_2} \)  
or there is no $(x(X), x^{m_1 m_2}(X))$-transverse intersection point for any hyperbolic metric~$X$.  
Since \( m_1 m_2 > 1 \), we have \( x \ne x^{m_1 m_2} \).  
It also follows that there is no $(x(X), y(X))$-transverse intersection point for any hyperbolic metric $X$.  
By Lemma~\ref{KeyLemma2}, we obtain \( i(x, y) = 0 \).

Assume that condition~(4) holds.  
We have
\begin{align}\label{condition4}
[x^m, c_1 y + c_2 y^{-1}] 
= c_1 [x^m, y] + c_2 [x^m, y^{-1}] = 0.
\end{align}
For a complete hyperbolic metric \( X \), let \( P \) be an $(x(X), y(X))$-transverse intersection point at which the forward angle is minimal. Then, we consider which term in equation~\eqref{condition4} may cancel the term \( |x^m_P y_P| \). There are two possibilities.

\textbf{Case (a)}  The cancelling term comes from the first bracket \( [x^m, y] \). 

\noindent
In this case, there exists another $(x(X), y(X))$-transverse intersection point \( Q \ne P \)  
such that \( |x^m_P y_P| = |x^m_Q y_Q| \).  
Then we have \( c_1(\varepsilon_P(x^m, y) + \varepsilon_Q(x^m, y)) = 0 \),  
and since \( c_1 \) is a non-zero divisor, it follows that \( \varepsilon_P(x^m, y) = -\varepsilon_Q(x^m, y) \).  
By Lemma~\ref{xy}, we obtain \( y = x^m \).  
In that case, equation~\eqref{condition4} becomes \( c_2 [x^m, x^{-m}] = 0 \),  
and since \( c_2 \) is also a non-zero divisor, we have \( [x^m, x^{-m}] = 0 \).  
Since \( x^m \ne x^{-m} \), Lemma~\ref{KeyLemma1} implies there is no $(x^{m}(X), x^{-m}(X))$-intersection point for any hyperbolic metric $X$. By this condition, it also follows that there is no $(x(X), y(X))$-transverse intersection point for any hyperbolic metric $X$. By Lemma~\ref{KeyLemma2}, we obtain $i(x, y) = 0.$

\textbf{Case (b)} The cancelling term comes from the second bracket \( [x^m, y^{-1}] \).

\noindent
In this case, there exists an $(x(X), y^{-1}(X))$-transverse intersection point \( Q \) such that \( |x^m_P y_P| = |x^m_Q y_Q^{-1}| \). Lemma~\ref{xyz} implies that there is no $(x(X), y(X))$-transverse intersection point for any hyperbolic metric $X$ since \( y \ne y^{-1} \). By Lemma~\ref{KeyLemma2}, we obtain $i(x, y) = 0.$ This completes the proof.
\end{proof}

\begin{remark}
Chas~\cite{Chas2004} constructed loops \( x, y \in \hat{\pi} \) satisfying \( [x, y] = 0 \), \( x \neq y \), and \( i(x, y) > 0 \).  
For any \( m \ge 2 \), using again the method introduced in the Appendix of Chas and Krongold~\cite{Chas-Krongold2010}, one can compute explicitly that \( [x^m, y] \neq 0 \).  
This observation is consistent with Theorem~\ref{SSC1} (Theorem \ref{SSC2}).
\end{remark}

\section{The center of the Goldman Lie algebra of a pair of pants}

In this section, we compute the center of the Goldman Lie algebra of a pair of pants. A key step is to understand when two loops of the form \( |x_P^m y_P| \) and \( |x_P^m z_P| \) can represent the same free homotopy class. The next lemma shows that, if this occurs for two distinct values of \(m\), then the lengths and forward angles of the original geodesics must coincide.

\begin{lemma} \label{Equality_lengths_angles}
Let $X$ be a complete hyperbolic metric, and $x, y, z \in \hat{\pi}$ be three pairwise distinct free homotopy classes of loops. Suppose that \(P\) is an \((x(X), y(X))\)-transverse intersection point and \(Q\) is an \((x(X), z(X))\)-transverse intersection point. If there exist distinct natural numbers \( m_1, m_2 \in \mathbb{N} \) such that
\[
|x_P^{m} y_P| = |x_Q^{m} z_Q| \quad \text{for } m = m_1, m_2,
\]
then we have
\[
\ell_y = \ell_z, \quad \phi_P = \phi_Q
\]
\end{lemma}

\begin{proof}
By Lemma \ref{thm_cosh}, for \( m = m_1\text{ and } m_2 \), we have
\begin{align*}
\coth\left( \frac{\ell_{x^m}}{2} \right) \cosh\left( \frac{\ell_y}{2} \right)
+ \sinh\left( \frac{\ell_y}{2} \right) \cos \phi_P
&= \coth\left( \frac{\ell_{x^m}}{2} \right) \cosh\left( \frac{\ell_z}{2} \right)
+ \sinh\left( \frac{\ell_z}{2} \right) \cos \phi_Q.
\end{align*}
Assume \( \ell_y \neq \ell_z \). Then we obtain
\[
\coth\left( \frac{\ell_{x^m}}{2} \right) =
\frac{
\sinh\left( \frac{\ell_z}{2} \right) \cos \phi_Q - \sinh\left( \frac{\ell_y}{2} \right) \cos \phi_P
}{
\cosh\left( \frac{\ell_y}{2} \right) - \cosh\left( \frac{\ell_z}{2} \right)
} \quad \text{for } m = m_1\text{ and } m_2.
\]
Since the right–hand side is independent of \(m\), this contradicts the strict monotonicity of the hyperbolic cotangent.
Therefore \(\ell_y = \ell_z\) and \(\phi_P = \phi_Q\).
\end{proof}

\begin{lemma} \label{CFLC}
Let $x, y_1, \ldots, y_k$ be pairwise distinct free homotopy classes of loops and $y = \sum_{j=1}^{k} c_j y_j$ where the coefficients $c_1, c_2, \ldots, c_k$ are nonzero elements of $K$. Then either $i(x, y_j) = 0$ for all $j \in \{1, \ldots, k\}$ or there exists a positive integer $m_0$ such that $[x^{m_0}, y] \neq 0$.
\end{lemma}

\begin{proof}
Assume that \( I \coloneqq i(x, y_1) + \cdots + i(x, y_k) \geq 1 \) and that \( [x^m, y] = 0 \) for all positive integers \( m \). By the definition of the Goldman bracket, for any positive integer $m \in \mathbb{N}$,
\[
[x^m, y] = \sum_{j=1}^{k} c_j [x^m, y_j] = m \sum_{j=1}^{k} c_j \sum_{P \in x \cap y_j} \varepsilon_{P}(x,y_j)|{x^m_{P}y_{j_P}}| .
\]
Here, the representatives of \( x^m \) and each \( y_j \) are taken to be the geodesic representatives with respect to a fixed complete hyperbolic metric \( X \). For each $m$, the sum
\[
\sum_{j=1}^k c_j \sum_{P \in x \cap y_j} \varepsilon_{P}(x,y_j)|{x^m_{P}y_{j_P}}|
\]
contains $I$ terms before cancellation. We may also assume $i(x, y_1) > 0$. Choose an \((x(X), y_1(X))\)-transverse intersection point \(P\) whose forward angle \(\phi_P\) is minimum among all \((x(X), y_1(X))\)-transverse intersection points,  
and consider the corresponding term \(|x^m_{P}y_{1_P}|\).  
In particular, since we have \( [x^m, y] = 0 \) for all \(1 \leq m \leq I + 1\),  
there exist \(1 \leq m_1 < m_2 \leq I+1\), \(y_j \in \hat{\pi}\),  
and an \((x(X), y_j(X))\)-transverse intersection point \(Q\) such that
\[
|x^m_{P}y_{1_P}| = |x^m_{Q}y_{j_Q}| \quad \text{for } m = m_1 \text{ and } m_2.
\]
By Lemma \ref{Equality_lengths_angles} and Lemma \ref{xyz}, there exists an  \((x(X), y_1(X))\)-transverse intersection point $R$ such that the forward angle at $R$ is strictly smaller than that at $P$, i.e., $\phi_R < \phi_P$. This contradicts the assumption that the forward angle $\phi_P$ was minimum among all \((x(X), y_1(X))\)-transverse intersection points. This completes the proof.
\end{proof}

The following lemma is classical.

\begin{lemma}\label{FU}
Let \(\Sigma\) be an orientable surface. If \(y \in \hat{\pi}\) is a free homotopy class of a loop on \(\Sigma\) such that \(i(x,y)=0\) for every \(x \in \hat{\pi}\), then \(y\) is non-essential.
\end{lemma}
 
In \cite[Theorem~7.2]{Kabiraj2016}, only simple closed curves $x$ are studied in view of Lemma~\ref{FU}. As a result, the statement fails when \( \Sigma \) is a pair of pants. By allowing arbitrary loops $x$, we can treat a pair of pants. For example, when \( \Sigma \) is a pair of pants, one can take $x$ to be a figure-eight curve.

\begin{thm}[Theorem \ref{COPOP1}] \label{COPOP2}
The center of the Goldman Lie algebra of a pair of pants is generated by the class of non-essential loops as a $K$-module.
\end{thm}

\begin{proof}
Let \( y = \sum_{j=1}^{k} c_j y_j \) be an element in the center, where the \( y_j \) are pairwise distinct free homotopy classes and each coefficient \( c_j \in K \) is nonzero. Let \( x \) be a figure-eight curve on the pair of pants. We fix an arbitrary direction of \( x \) to compute the Goldman bracket. Since \( y \) is central, we have \([x^m, y] = 0\) for all \( m \in \mathbb{N} \). By Lemma~\ref{CFLC} and Lemma~\ref{FU}, each \( y_j \) satisfies $i(x, y_j) = 0$ and so is non-essential. This proves the theorem.
\end{proof}

\begin{remark}
In \cite{Kabiraj2023}, Kabiraj proved that the Poisson centers of the symmetric algebra \( S(K\hat{\pi}) \) and the universal enveloping algebra \( U(K\hat{\pi}) \) are generated by the classes of non-essential loops, for any surface (with or without boundary) except a pair of pants. Here, \( S(K\hat{\pi}) \) denotes the symmetric algebra and \( U(K\hat{\pi}) \) denotes the universal enveloping algebra, both of which naturally carry Poisson algebra structures. Thanks to Lemma~\ref{xyz}, a similar proof applies to the case of a pair of pants, and we can conclude that the Poisson centers of \( S(K\hat{\pi}) \) and \( U(K\hat{\pi}) \) are also generated by the classes of non-essential loops in this case. Moreover, in Sections~7--9 of \cite{Kabiraj2023}, Kabiraj also proved that the Poisson centers of various types of skein algebras are generated by the classes of non-essential loops, for any surface (with or without boundary) except a pair of pants. These results follow immediately once the Poisson centers of the symmetric algebra \( S(K\hat{\pi}) \) and the universal enveloping algebra \( U(K\hat{\pi}) \) are determined. Therefore, as we have established the corresponding result for a pair of pants, it follows that the Poisson centers of various skein algebras for a pair of pants are also generated by the classes of non-essential loops.
\end{remark}

\bibliography{SC}
\bibliographystyle{plain}

\end{document}